\documentclass[11pt]{amsart}
\usepackage{amsmath, amssymb, comment, epsfig}
\usepackage{psfrag}

\newcommand{\K}{\mathcal{K}}
\newcommand{\Gr}{\mathcal{G}r}
\newcommand{\wl}{X^*}
\newcommand{\cwl}{X_*}
\newcommand{\D}{\mathrm{D}}
\newcommand{\C}{\mathbb{C}}
\newcommand{\Z}{\mathbb{Z}}
\newcommand{\R}{\mathbb{R}}
\newcommand{\N}{\mathbb{N}}
\newcommand{\val}{\operatorname{val}}

\newcommand{\row}{\Gamma}
\newcommand{\conv}{\operatorname{conv}}
\newcommand{\bigdot}{\bullet}
\newcommand{\wi}{\mathbf{i}}

\newcommand{\realt}{\mathfrak{t}_\R}
\newcommand{\fund}{\Lambda}
\newcommand{\Hom}{\operatorname{Hom}}

\newcommand{\Lieg}{\mathfrak{g}}
\newcommand{\wt}{\operatorname{wt}}
\newcommand{\Jl}{v}
\newcommand{\Comp}{\operatorname{Comp}}

\newtheorem{Theorem}{Theorem}[section]
\newtheorem{Proposition}[Theorem]{Proposition}
\newtheorem{Lemma}[Theorem]{Lemma}
\newtheorem{Corollary}[Theorem]{Corollary}
\newtheorem{Conjecture}[Theorem]{Conjecture}
\newtheorem{Question}{Question}
\newtheorem{Specialthm}{Theorem}

\theoremstyle{definition}
\newtheorem{Example}[Theorem]{Example}

\author{Joel Kamnitzer}
\title{The crystal structure on the set of Mirkovi\'c-Vilonen Polytopes}
\address{American Institute of Mathematics \\ 360 Portage Ave \\ Palo Alto, CA 94306}
\email{jkamnitz@aimath.org}
\date{\today}

\begin{document}

\begin{abstract}
In an earlier work, we proved that MV polytopes parameterize both Lusztig's canonical basis and the Mirkovi\'c-Vilonen cycles on the Affine Grassmannian.  Each of these sets has a crystal structure (due to Kashiwara-Lusztig on the canonical basis side and due to Braverman-Finkelberg-Gaitsgory on the MV cycles side).  We show that these two crystal structures agree.  As an application, we consider a conjecture of Anderson-Mirkovi\'c which describes the BFG crystal structure on the level of MV polytopes.  We prove their conjecture for $ \mathfrak{sl}_n $ and give a counterexample for $\mathfrak{sp}_6$.  Finally we explain how Kashiwara data can be recovered from MV polytopes.
\end{abstract}

\maketitle
\tableofcontents

\section{Introduction}

Let $ \Lieg^\vee $ be a complex semisimple Lie algebra.  A $ \Lieg^\vee$-crystal is a combinatorial object corresponding to a representation of $ \Lieg^\vee $.  A particularly important crystal is the crystal $ B(\infty) $ corresponding to the Verma module of $ \Lieg^\vee $ of highest weight $ 0 $.  There are a number of ways to realize this crystal.  Historically, the first construction of this crystal used Kashiwara's crystal basis for $ U^\vee_- $ as the underlying set of the crystal.  This description is representation theoretic, involving the upper triangular part of the quantized universal enveloping algebra of $ \mathfrak{g}^\vee $.  

More recently, Braverman-Finkelberg-Gaitsgory \cite{BG,BFG} gave a geometric description of this crystal where the underlying set is the set of Mirkovi\'c-Vilonen cycles on the affine Grassmannian for the group $ G $, where $ G $ is the simply connected, semisimple group dual to $ \Lieg^\vee $.  In \cite{jared2}, Anderson proposed studying these cycles by means of their moment map images, which he called MV polytopes.  

In \cite{me}, we gave an explicit characterization of these polytopes and showed that they describe both MV cycles and Lusztig's canonical basis for $ U^\vee_- $ (which is the same as Kashiwara's crystal basis).  Thus, we have bijections
\begin{equation*} \label{eq:bij}
\mathcal{B} \longleftrightarrow \mathcal{P} \longleftrightarrow \mathcal{M}
\end{equation*}
where $ \mathcal{B} $ denotes the canonical basis, $ \mathcal{P} $ denotes the set of MV polytopes, and $ \mathcal{M} $ denotes the set of MV cycles.

Our main result is the following.

\begin{Specialthm}[Theorem \ref{th:LBZ=BFG}]
Via these bijections, we get two crystal structures on the set of MV polytopes.  These two crystal structures agree.
\end{Specialthm}

We believe that this theorem shows the naturality and importance of the crystal structure on the set of MV polytopes.  The idea of constructing $ B(\infty) $ using the set of MV polytopes  is due to Anderson-Mirkovi\'c.

Alternatively, one can say Theorem \ref{th:LBZ=BFG} proves that the above bijection $ \mathcal{B} \rightarrow \mathcal{M} $ is the unique isomorphism of crystals between the canonical basis and the set of MV cycles (the existence and uniqueness of such an isomorphism follows from a uniqueness theorem for $B(\infty) $, see \cite{BFG,KS}).  The explicit construction of this isomorphism answers a question posed by Braverman-Gaitsgory \cite{BG}.

\begin{figure}
\begin{center}
\psfrag{a2}{$\scriptstyle{-\alpha_1^\vee}$} \psfrag{a1}{$\scriptstyle{-\alpha_2^\vee}$}
\epsfig{file=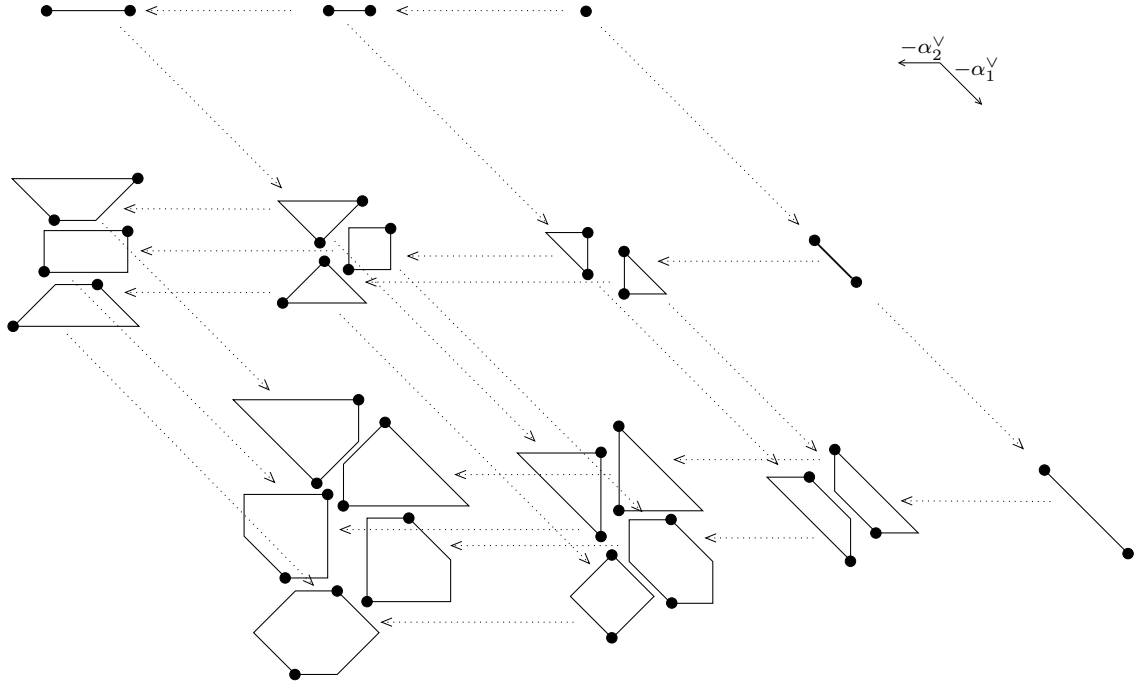, height=9cm}
\caption{A portion of the $ \mathfrak{sp}_4 $ $B(\infty)$ crystal.}
\label{fig:crys}
\end{center}
\end{figure}

The crystal structure on MV polytopes is fairly easy to describe. Let $ P $ be an MV polytope.  Since an MV polytope is a pseudo-Weyl polytope, it comes with a map $ w \mapsto \mu_w $ from the Weyl group onto its vertices (see \ref{se:pWeyl}).  The crystal operator $ f_j $ acts on $ P $ to produce the unique MV polytope with the vertex $ \mu_w $ unchanged if $ s_j w < w $ and with the vertex $ \mu_e $ shifted by $ - \alpha_j^\vee $.  The rest of the vertices of $f_j \cdot P $ are determined by the tropical Pl\"ucker relations (section \ref{se:MVpoly}).  Figure \ref{fig:crys} shows with a portion of the crystal graph of MV polytopes for $ G = Sp_4$.

This description is immediate from the canonical basis side (Theorem \ref{th:LBZ}).  We simply combine the way that Lusztig data of canonical basis elements changes under crystal operators (Proposition \ref{th:cryscan}) with the relationship between the Lusztig datum of a canonical basis element and the lengths of edges of the corresponding MV polytopes (Theorem \ref{th:can}). 

It is more difficult to show that the MV cycle world gives rise to the same crystal structure.  The crystal operator $ f_j $ of \cite{BFG} is defined by factoring an open part of each MV cycle $Z $ as a product $ A \times B $, where $ A$ is an MV cycle for the rank 1 Levi subgroup $ G_j $ and $ B $ is a ``relative MV cycle''.  The action of $ f_j $ is defined by its action on $ A $ (sections \ref{se:restrict}, \ref{se:BFGcrys}).  To relate this geometric procedure to the combinatorial one outlined above, we show that the generic values of certain constructible functions on $ Z $ only depend on their generic values on $ B $ (Proposition \ref{th:relval}) and hence are invariant under the crystal operator.  On the MV polytope level, this shows that certain vertices (namely the ones described above) do not change under the crystal operator. This allows us to establish the equivalence of the two crystal structures (Theorem \ref{th:LBZ=BFG}).

\subsection{Anderson-Mirkovic conjecture}

Inspired by the BFG crystal structure on MV cycles, Anderson-Mirkovi\'c conjectured a crystal structure on MV polytopes (see section \ref{se:AM}).  This crystal structure is similar, but more specific than the one described above (see Proposition \ref{th:AMBZ}, Example \ref{eg:AM}, and Corollary \ref{th:AMtypeA}), so it is very interesting to examine the validity of this conjecture.  Their conjecture was one of the main motivations for this work.

\begin{Specialthm}
The Anderson-Mirkovic conjecture holds for $ \mathfrak{sl}_n$ but fails for $ \mathfrak{sp}_6 $.
\end{Specialthm}

To prove that the conjecture holds for $ \mathfrak{sl}_n $ we use a certain inductive argument which allows us to use the tropical Plucker relations in a systematic manner (Theorem \ref{th:jclosetrue}).  For $ \mathfrak{sp}_6 $ we present a counterexample in section \ref{se:counter}.

\subsection{Additional combinatorial structure}

The crystal $ B(\infty) $ has additional combinatorial structure.  We close the paper by considering how some of this combinatorial structure is reflected in the MV polytope model.

Kashiwara \cite{Kas} showed that there is an involution $ * $ of the set $ B(\infty) $ and hence an additional family of crystal operators defined by twisting the original crystal operators by the Kashiwara involution.  We show that in the MV polytope model, the Kashiwara involution corresponds to the map $ P \mapsto -P $ of negating a polytope and that the twisted crystal structure corresponds to acting on the highest vertex as opposed to the lowest vertex.

Another combinatorial structure is the notion of Kashiwara data.  Consider the crystal $B(\lambda) $ corresponding to the finite dimensional representation of highest weight $ \lambda $. The Kashiwara (or string) datum of an element of this crystal with respect to a reduced word $ \mathbf{i} $ for $ w_0 $ is the sequence of lengths of steps needed to reach the bottom of the crystal by applying the lowering operators $ f_{i_1}, f_{i_2}, \dots $.  This notion was introduced by Kashiwara \cite{Kas} and Berenstein-Zelevinsky \cite{BZstring}.

The MV polytope model for $ B(\lambda) $ consists of those polytopes whose highest weight is $ \lambda $ and which fit inside $ \conv (W \cdot \lambda) $.   These MV polytopes encode their own Kashiwara data in a natural way.  
\begin{Specialthm}
A reduced word $ \mathbf{i} $ determines a path along the edge of any MV polytope.
The Kashiwara datum of an MV polytope $ P \in B(\lambda) $ with respect to $ \mathbf{i} $ is the sequence of heights of the midpoints of the edges along this path.
\end{Specialthm}
This result was inspired by Morier-Genoud's theorem \cite{MG} describing the Sch\"utzenberger involution in terms of Lusztig and Kashiwara data.  However, we do not make use of her result.  Our proof uses Berenstein-Zelevinsky's observation \cite{BZtpm} relating change of parametrization of Kashiwara data with the tropical Pl\"ucker relations.

Note that there is another family of Kashiwara data given by considering raising crystal operators.  In a sense, this is the more natural family to consider since it gives the notion of Kashiwara data on $ B(\infty) $.  Unfortunately, this raising Kashiwara data is not easily expressed in terms of MV polytopes --- in particular it is not a linear function of the BZ datum $ M_\bigdot $.  To describe this Kashiwara data one needs to use double BZ data indexed by the set $ \bigcup_i W \cdot \fund_i \times W \cdot \fund_i $ as explained in \cite{BZtpm}.  We do not know of any way to relate this double BZ data to MV polytopes.

\subsection*{Acknowledgements}
The encouragement and help of my advisor Allen Knutson was important at many stages of this project.  I thank Jared Anderson for explaining his conjecture to me and for useful conversations about MV polytopes.  I also benefited from conversations with Arkady Berenstein, Alexander Braverman, Misha Kogan, Alexander Postnikov, and Peter Tingley.  I am also grateful to Allen Knutson and Peter Tingley for their careful reading of this text.  In the course of this work, I was supported by an NSERC postgraduate scholarship and an AIM five-year fellowship.

\section{Notation}
Let $ G $ be a connected simply-connected semisimple complex group. 

Let $ T $ be a maximal torus of $ G $ and let $ \wl = \Hom(T, \C^\times), \cwl = \Hom(\C^\times, T) $ denote the weight and coweight lattices of $ T $.  Let $ \Delta \subset \wl $ denote the set of roots of $ G $.  Let $ W = N(T)/ T $ denote the Weyl group.

Let $ B $ be a Borel subgroup of $ G $ containing $ T $.  Let $ \alpha_1, \dots, \alpha_r $ and $ \alpha_1^\vee, \dots, \alpha_r^\vee $ denote the simple roots and coroots of $ G $ with respect to $ B $.  Let $ N $ denote the unipotent radical of $ B $.  Let $ \fund_1, \dots, \fund_r $ be the fundamental weights.  Let $ I = \{ 1, \dots, r \} $ denote the vertices of the Dynkin diagram of $ G $.  Let $ a_{ij} = \langle \alpha_i^\vee, \alpha_j \rangle $ denote the Cartan matrix.  

Let $ s_1, \dots, s_r \in W $ denote the simple reflections.  Let $ e $ denote the identity in $ W $ and let $ w_0 $ denote the longest element of $ W $.  Let $ m $ denote the length of $ w_0 $ or equivalently the number of positive roots.  We will also need the Bruhat order on $ W $, which we denote by $ \ge $. 

We also use $ \ge $ for the usual partial order on $ \cwl $, so that $ \mu \ge \nu $ if and only if $ \mu - \nu $ is a sum of positive coroots.  More generally, we have the twisted partial order $ \ge_w $, where $ \mu \ge_w \nu $ if and only if $ w^{-1}  \cdot \mu \ge w^{-1} \cdot \nu $.

Let $ \realt := \cwl \otimes \R $ (the Lie algebra of the compact form of $ T $).  For each $ w$, we extend $ \ge_w $ to a partial order on $ \realt $, so that $ \beta \ge_w \alpha $ if and only if $ \langle \beta - \alpha, w \cdot \fund_i \rangle \ge 0 $ for all $ i $.

For each $ i \in I $, let $ \psi_i : SL_2 \rightarrow G $ be denote the $ i$th root subgroup of $ G $.

For $ w \in W $, let $ \overline{w} $ denote the lift of $ w $ to $ G $, defined using the lift of $ \overline{s_i} := \psi_i \Big( \big[ \begin{smallmatrix} 0 & 1 \\ -1 & 0 \\ \end{smallmatrix} \big] \Big) $.

A \textbf{reduced word} for an element $ w \in W $ is a sequence of indices $ \textbf{i} = (i_1, \dots, i_k) \in I^k$ such that $ w = s_{i_1} \cdots s_{i_k} $ is a reduced expression.  In this paper, reduced word will always mean reduced word for $ w_0 $.

If $ X $ is any variety, we write $ \Comp(X) $ for the set of components of $ X $.

\section{Crystal structure via canonical basis}

We begin by reviewing some notions from \cite{me} concerning MV polytopes.

\subsection{Pseudo-Weyl polytopes} \label{se:pWeyl}
A Weyl group translate $ w \cdot \fund_i $ of a fundamental weight is called a \textbf{chamber weight} of level $ i $.  The collection of all chamber weights is denoted $ \Gamma = \{ w \cdot \fund_i : w \in W, i \in I \} $.  

We say that a collection $ M_\bigdot = \big( M_\gamma \big)_{ \gamma\in \row } $ of integers satisfies the \textbf{edge inequalities} if
\begin{equation} \label{eq:nondeg}
M_{w \cdot \fund_i} + M_{w s_i \cdot \fund_i} + \sum_j a_{ji} M_{w \cdot \fund_j} \le 0, 
\end{equation} 
for all $ i \in I $ and $ w\in W $.

\begin{figure}
\begin{center}
\psfrag{me}{$\mu_e$} \psfrag{ms1}{$\mu_{s_1}$} \psfrag{ms1s2}{$\mu_{s_1 s_2}$}
\psfrag{ms1s2s1}{$\mu_{s_1s_2s_1}$} \psfrag{ms2s1}{$\mu_{s_2s_1}$} \psfrag{ms2}{$\mu_{s_2}$}
\psfrag{f1}{$\scriptstyle{\fund_1}$} \psfrag{s1f1}{$\scriptstyle{s_1 \cdot \fund_1}$} \psfrag{s1s2f2}{$\scriptstyle{s_1 s_2 \cdot \fund_2}$}
\psfrag{f2}{$\scriptstyle{\fund_2}$} \psfrag{s2f2}{$\scriptstyle{s_2 \cdot \fund_2}$} \psfrag{s2s1f1}{$\scriptstyle{s_2 s_1 \cdot \fund_1}$}
\epsfig{file=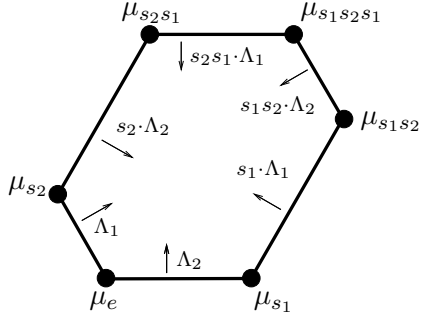, height=4cm}
\caption{A pseudo-Weyl polytope for $SL_3 $.}
\label{fig:pWeylsl3}
\end{center}
\end{figure}

Given such a collection, we can form the \textbf{pseudo-Weyl polytope} $ P(M_\bigdot) := \{ \alpha \in \realt : \langle \alpha, \gamma \rangle \ge M_\gamma $ for all $ \gamma \} $.  In \cite{me}, we showed that such polytopes come with a map $ w \mapsto \mu_w $ from the Weyl group onto the vertices of the polytope such that 
\begin{equation} \label{eq:mufromM}
\langle \mu_w, w \cdot \fund_i \rangle = M_{w \cdot \fund_i} .
\end{equation}
This collection of coweights $ \mu_\bigdot = \big( \mu_w \big)_{w \in W} $ is called the \textbf{GGMS datum} of the pseudo-Weyl polytope.  Moreover, we have that 
\begin{equation*}
P(M_\bigdot) = \{ \alpha \in \realt : \alpha \ge_w \mu_w \text{ for all } w \in W \}.
\end{equation*}
Figure \ref{fig:pWeylsl3} shows an example of a pseudo-Weyl polytope for $ G = SL_3 $ with vertices and chamber weights labelled.

\subsection{MV polytopes} \label{se:MVpoly}
Let $ w \in W, i,j \in I $ be such that $ w s_i > w, ws_j > w $, and $ i \ne j$.  We say that a collection $ \big( M_\gamma \big)_{\gamma \in \Gamma} $ satisfies the \textbf{tropical Pl\"ucker relation} at $ (w, i,j) $ if $a_{ij} = 0 $ or if $ a_{ij} = a_{ji} = -1 $ and
\begin{equation} \label{eq:A2trop}
M_{ws_i \cdot \fund_i} + M_{w s_j \cdot \fund_j} = \min(M_{w \cdot \fund_i} + M_{w s_i s_j \cdot \fund_j} , M_{w s_j s_i \cdot \fund_i} + M_{w \cdot \fund_j} ) .  
\end{equation}
There are also conditions for the other possible values of $ a_{ij}, a_{ji} $.  We leave out these cases since they will not be used in this paper.  See \cite{me} for these cases.

We say that a collection $ M_\bigdot = \big( M_\gamma \big)_{\gamma \in \Gamma} $ satisfies the \textbf{tropical Pl\"ucker relations} if it satisfies the tropical Pl\"ucker relation at each $ (w, i, j) $. 

A collection $ M_\bigdot $ of integers is called a \textbf{BZ datum} of coweight $ (\mu_1, \mu_2) $ if:
\begin{enumerate}
\item $ M_\bigdot $ satisfies the tropical Pl\"ucker relations.
\item $ M_\bigdot $ satisfies the edge inequalities (\ref{eq:nondeg}).
\item If $ \mu_\bigdot $ is the GGMS datum of $ P(M_\bigdot) $, then $ \mu_e = \mu_1 $ and $\mu_{w_0} = \mu_2 $.
\end{enumerate}

If $ M_\bigdot $ is a BZ datum of coweight $(\mu_1, \mu_2) $, then $P(M_\bigdot) $ is called an \textbf{MV polytope} of coweight $(\mu_1, \mu_2$).  

In \cite{jared2}, Anderson defined MV polytopes to be the moment map images of the MV cycles.  In \cite{me}, using that definition, we proved that they were of the form $ P(M_\bigdot) $ for $ M_\bigdot $ a BZ datum.  Here, we are starting with the polytopes, so we prefer to define them this way.

There is an action of $ \cwl $ on $\realt$ by translation.  If $ M_\bigdot $ is a BZ datum of coweight $(\mu_1, \mu_2) $ then $ \nu + P(M_\bigdot) = P(M'_\bigdot) $ where $ M'_\gamma = M_\gamma + \langle \nu, \gamma \rangle $.  In particular, $ M'_\bigdot $ is a BZ datum of coweight $ (\mu_1 + \nu, \mu_2 + \nu) $. Hence we have an action of $ \cwl $ on the set of BZ datum and on the set of MV polytopes.  

The orbit of an MV polytope of coweight $ (\mu_1, \mu_2) $ is called a \textbf{stable MV polytope} of coweight $ \mu_1 - \mu_2 $.  Note that a stable MV polytope of coweight $\mu $ has a unique representative of coweight $ (\nu + \mu, \nu) $ for all coweights $ \nu $.  Let $ \mathcal{P} $ denote the set of stable MV polytopes.

Let $ \wi = (i_1, \dots, i_m) $ be a reduced word for $ w_0 $.  Let $w^\wi_k = s_{i_1} \cdots s_{i_k} $.  We define the $\wi$-\textbf{Lusztig datum} of any MV polytope $ P(M_\bigdot) $ to be the vector of non-negative integers $ n_\bigdot = (n_1, \dots, n_m) $, defined by
\begin{equation} 
n_k = - M_{w_{k-1}^\wi \cdot \fund_{i_k}} - M_{w_k^\wi \cdot \fund_{i_k}} - \sum_{j \ne i} a_{ji} M_{w_k^\wi \cdot \fund_j}. 
\end{equation}
As explained in \cite{me}, $\wi $ determines a path $ e = \mu_{w^\wi_0}, \mu_{w^\wi_1}, \dots, \mu_{w_m^\wi} = w_0 $ through the 1-skeleton of the polytope.  The integers $ (n_1, \cdots, n_m) $ are the lengths of the edges along this path, in the sense that 
\begin{equation} \label{eq:lengths}
\mu_{w_k^\wi} - \mu_{w_{k-1}^\wi} = n_k w^\wi_{k-1} \cdot \alpha_{i_k}^\vee. 
\end{equation}

Note that the $ \wi$-Lusztig datum of an MV polytope is invariant under the action of $\cwl$ and hence the $ \wi$-Lusztig datum of a stable MV polytope is well-defined.  On the other hand, the GGMS datum $ \mu_\bigdot $ of a stable MV polytope is defined only up to simultaneous shift $ \mu_\bigdot \mapsto \mu_\bigdot + \nu $.  

\begin{Proposition}[\cite{me}] \label{th:lusdatbij}
For any reduced word $ \wi $, the map $ \mathcal{P} \rightarrow \N^m $ defined by taking $ \wi$-Lusztig datum is a bijection.
\end{Proposition}

\subsection{Canonical basis}

Recall that $ \mathfrak{g}^\vee $ is the Lie algebra whose root system is dual to that of $ G $.  In particular, the weight lattice of $ \mathfrak{g}^\vee $ is $ \cwl $.  Let $ \mathcal{B} $ denote Lusztig's canonical basis for $ U^\vee_-$, the lower triangular part of the quantized universal enveloping algebra of $ \mathfrak{g}^\vee $.  Lusztig showed that a choice of reduced word $ \wi $ for $ w_0 $ gives rise to a bijection $ \phi_{\wi} : \mathcal{B} \rightarrow \N^m $ (see \cite[section 2]{Lbook} or \cite[Proposition 4.2]{BZtpm} for more details). Following Berenstein-Zelevinsky, we call $\phi_{\wi}(b) $ the $\wi$-\textbf{Lusztig datum} of $ b $.
  
In \cite{me}, we proved that MV polytopes parameterize the canonical basis.  This used the work of Berenstein-Zelevinsky \cite{BZtpm}.  
  
\begin{Theorem} \label{th:can}
There is a coweight preserving bijection $ b \mapsto P(b) $ between the canonical basis $ \mathcal{B} $ and the set $\mathcal{P} $of stable MV polytopes.  Under this bijection, the $\wi$-Lusztig datum of $ b $ equals the $\wi$-Lusztig datum of $P(b)$.
\end{Theorem}
In other words to find the $\wi$-Lusztig datum of $ b $, we can just look at the lengths of the edges in $P(b) $ along the path determined by $ \wi $.

\subsection{Crystal structure on the canonical basis}
For our purposes, a \textbf{crystal} is a set $ \mathcal{C} $ along with structure maps
\begin{equation*}
e_j : \mathcal{C} \rightarrow \mathcal{C} \cup \{ 0 \}, \ f_j : \mathcal{C} \rightarrow \mathcal{C} \cup \{ 0 \}, \text{ for each $ j \in I $, and } \wt : \mathcal{C} \rightarrow \cwl,
\end{equation*}
which satisfy the following axioms.
\begin{enumerate}
\item If $ b \in \mathcal{C} $, $ j \in I $, and $ e_j \cdot b \ne 0 $, then $\wt(e_j \cdot b) = \wt(b) + \alpha^\vee_j $. 
\item If $ b \in \mathcal{C} $, $ j \in I $, and $ f_j \cdot b \ne 0 $, then $\wt(f_j \cdot b) = \wt(b) - \alpha^\vee_j$.
\item $ b' = e_j \cdot b $ if and only if $ f_j \cdot b' = b $.
\end{enumerate}

The ``weight'' function $ \wt $ takes values in $\cwl$, the weight lattice of $ \mathfrak{g}^\vee$, so that these are actually crystal for representations of $ \mathfrak{g}^\vee $.

Note that if we are dealing with the crystal $ B(\infty) $, then we will never have $ f_j \cdot b = 0 $.

\subsection{Crystal structure on the canonical basis}
Following from Kashiwara's work on crystal bases, Lusztig proved the following result.
\begin{Theorem}[\cite{Lbookbook}]
There exists a crystal whose underlying set is Lusztig's canonical basis $ \mathcal{B}$ and whose $\wt $ function coincides with the usual coweight function on the canonical basis.
\end{Theorem}

For our purposes, the key aspect of this crystal structure is that, as long as $\wi$ is chosen correctly, it is easy to see how the $\wi$-Lusztig datum changes under the application of $ f_j $.
\begin{Proposition}[\cite{Lbookbook}, {\cite[Prop. 3.6]{BZtpm}}] \label{th:cryscan}
Let $ j \in I$ and let $ \wi $ be a reduced word with $i_1 = j $.  Let $ b \in \mathcal{B} $ and let $ n_\bigdot $ be its $\wi$-Lusztig datum.  Then,
\begin{enumerate}
\item $ f_j \cdot b $ has $\wi$-Lusztig datum $(n_1 + 1,n_2 \dots, n_m) $.
\item If $n_1 = 0 $, then $ e_j \cdot b = 0 $.  Otherwise, $e_j \cdot b $ has $\wi $-Lusztig datum $(n_1 - 1, n_2,\dots, n_m)$.
\end{enumerate}
\end{Proposition}

\subsection{Crystal structure on MV polytopes}
Since we have a bijection (Theorem \ref{th:can}) between the canonical basis and the set of stable MV polytopes, the crystal structure on the canonical basis gives a crystal structure on $ \mathcal{P}$.  We call this the LBZ (Lusztig-Berenstein-Zelevinsky) crystal structure.

To describe this crystal structure, we will choose a representative $ [P] \in \mathcal{P} $ of some coweight $ (\mu_1, \mu_2) $.  Then we will define $ f_j \cdot P $ and $ e_j \cdot P $ as MV polytopes of coweights $ (\mu_1 - \alpha_j^\vee, \mu_2) $ and $(\mu_1 + \alpha_j^\vee, \mu_2 ) $.  In other words, we will always be moving the bottom vertex $ \mu_e $ of the MV polytope.  Note that $ [f_j \cdot P], [e_j \cdot P] $ will be stable MV polytopes of coweights $ \mu_1 - \mu_2 \mp \alpha_j^\vee $ as expected. 

\begin{Theorem} \label{th:LBZ}
Let $ P $ be an MV polytope with GGMS datum $ \mu_\bigdot $.
\begin{enumerate}
\item $ f_j \cdot P $ is the unique MV polytope whose set of vertices $ \mu'_\bigdot $ satisfies
\begin{equation*}
\mu'_e = \mu_e - \alpha_j^\vee \ \text{ and } \ \mu'_w = \mu_w \text{ if } s_j w < w.
\end{equation*}
\item $ e_j \cdot P = 0 $ if and only if $ \mu_e = \mu_{s_j} $.  Otherwise, $ e_j \cdot P $ is the unique MV polytope whose set of vertices $ \mu'_\bigdot $ satisfies
\begin{equation*}
\mu'_e = \mu_e + \alpha_j^\vee \ \text{ and } \ \mu'_w = \mu_w \text{ if } s_j w < w.
\end{equation*}
\end{enumerate}
\end{Theorem}

\begin{proof}
Note that it suffices to prove the half of this theorem that deals with $ f_j$.  Also note that the MV polytope $ f_j \cdot P $ is determined by its vertices $ \{ \mu'_w : s_j w < w \} \cup \{ \mu'_e \} $ since to determine an MV polytope we only need to give the vertices along a path corresponding to a particular reduced word (this follows from Theorem \ref{th:lusdatbij}).  Moreover, by assumption, $ \mu'_e = \mu_e - \alpha_j^\vee $.  Hence it suffices to prove that for all $ w \in W $ with $ s_j w < w $, then $ \mu'_w = \mu_w $.

Let $ w \in W $ be such that $ s_j w < w $.  Then there exists a reduced word $\wi$ and an $ l $ such that $ \wi_1 = j $ and $ w = w_l^\wi $.  Let $ (n_1, \dots, n_m) $ be the $\wi$-Lusztig datum of $ P $.  By Proposition \ref{th:cryscan}, the $\wi$-Lusztig datum of $ f_j \cdot P $ is $(n_1 + 1, n_2, \dots, n_m)$.  Then the length formula (\ref{eq:lengths}), shows that $ \mu'_{w_k^\wi} = \mu_{w_k^\wi} $ for $ k = 1, \dots, l $.  In particular, $ \mu'_w = \mu_w $ as desired.
\end{proof}

Theorem \ref{th:LBZ} implies that if $ P = P(M_\bigdot) $ and if $ f_j \cdot P = P(M'_\bigdot) $, then $ M'_{\fund_j} = M_{\fund_j} - 1 $ and $ M'_\gamma = M_\gamma  $ whenever $ \gamma = w \cdot \fund_i $ for some $ w $ such that  $ s_j w < w $ (these $\gamma $ are called $ j$-\textbf{relative} chamber weights).  The rest of the $ M'_\gamma $ (and hence the rest of the vertices of $ f_j \cdot P $) are determined by the tropical Pl\"ucker relations.  However, this is quite non-explicit as it involves recursively solving many $(\min, +)$ equations.  Hence it is interesting to have a formula which explicitly determines the rest of the $ M'_\gamma $.  The Anderson-Mirkovi\'c conjecture which we examine in section \ref{se:AM} can be viewed a solution to this problem (also see Proposition \ref{th:AMBZ}).  In particular, for $ \Lieg^\vee = \mathfrak{sl}_n $, we give an explicit formula in Corollary \ref{th:AMtypeA}.

\begin{Example} \label{eg:LBZ}
Consider $ G = SL_3 $.  Then $ \cwl = \{ (a,b,c) \in \Z^3 : a + b + c = 0 \} $ and $ \wl = \Z^3 / \Z \cdot (1,1,1) $.  We can write $  \Gamma = \{ 1, 2, 3, 12, 13, 23 \} $ since we identify chamber weights with proper subsets of $\{1,2,3 \} $ (for example, $ 12 $ is the weight $ (1,1,0) \in \wl $).  Consider the MV polytope $ P $ with BZ datum
\begin{equation*}
M_1 = M_2 =  M_3 =  M_{12} =  M_{13} = M_{23} = -1.
\end{equation*}
This polytope has GGMS datum
\begin{gather*}
\mu_e = (-1,0,1), \ \mu_{s_1} =  (0, -1, 1), \ \mu_{s_2} =  (-1, 1, 0), \\ \mu_{s_1 s_2} =  (1, -1, 0), \ \mu_{s_2 s_1} = (0, 1, -1),\  \mu_{w_0} =  (1, 0, -1).
\end{gather*} 

\begin{center}
\epsfig{file=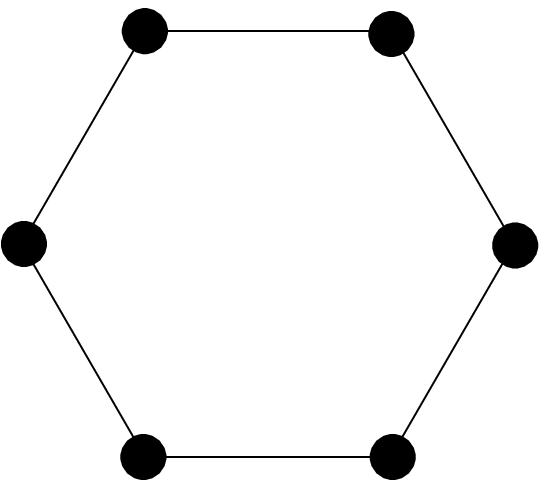,height=2cm}
\end{center}

So if $ M'_\bigdot, \mu'_\bigdot$, denote the BZ and GGMS data of $ f_1 \cdot P $, then by Theorem \ref{th:LBZ},
\begin{gather*}
\mu'_e = \mu_e - \alpha_1^\vee = (-2, 1, 1), \ \mu'_{s_1} = \mu_{s_1} = (0,-1, 1), \\ \mu'_{s_1 s_2} = \mu_{s_1 s_2} = (1,-1,0), \ \mu'_{w_0} = \mu_{w_0} = (1,0,-1), \\ 
M'_1= M_1 - 1 = -2, \ M'_2 = M_2 = -1, \ M'_3 = M_3 = -1, \\ M'_{12}= M_{12} = -1, \ M'_{23} = M_{23} = -1.
\end{gather*}

\begin{center}
\epsfig{file=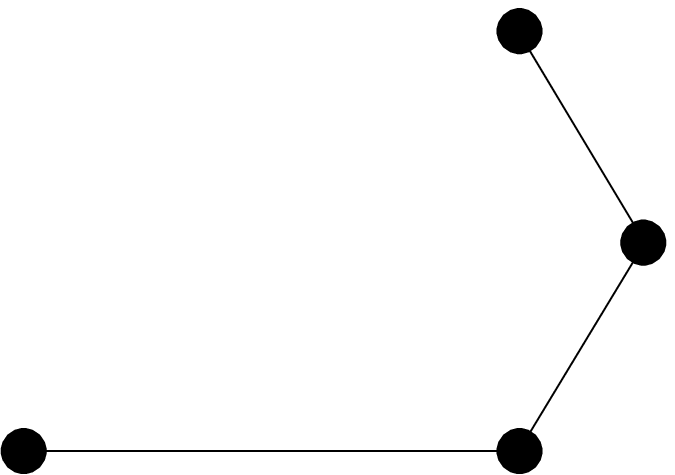,height=2cm}
\end{center}

Now, we use the tropical Pl\"ucker relation (see (\ref{eq:A2trop}))
\begin{equation*}
M'_{13} + M'_2 = \min( M'_1 + M'_{23}, M'_3 + M'_{12} ),
\end{equation*}
to find that $ M'_{13} = -2 $.  Hence $ \mu'_{s_2} = (-2,2, 0) $ and $ \mu'_{s_2 s_1} = (-1,2,-1) $.
\begin{center}
\epsfig{file=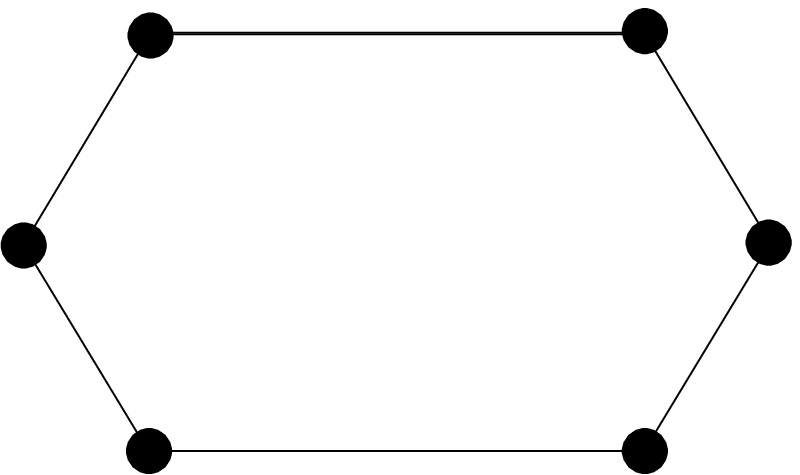,height=2cm}
\end{center}
\end{Example}

\section{Crystal structure via MV cycles}
\subsection{MV cycles}
We begin by recalling some definitions from \cite{me}, where we give more discussion about them.

Let $ \K = \C((t)) $ denote the field of Laurent series and let $ \mathcal{O} = \C[[t]] $ denote the ring of power series.  We define  the \textbf{affine Grassmannian} to be the left quotient $ \Gr = G(\mathcal{O}) \setminus G(\K) $.  Note that we have a right action of $ G(\K) $ on $\Gr $.

A coweight $ \mu \in \cwl $ gives a homomorphism $ \mathbb{C}^\times \rightarrow T $ and hence an element of $ \Gr$.  We denote the corresponding element $ t^\mu $. 

For $ w \in W $, let $ N_w = w N w^{-1} $.  For $ w \in W $ and $ \mu \in \cwl $ define the \textbf{semi-infinite cells}
\begin{equation}
 S_w^\mu := t^\mu N_w(\K). 
\end{equation}

Given any collection $ \mu_\bigdot = \big( \mu_w \big)_{w \in W} $ of coweights, we can form the \textbf{GGMS stratum}
\begin{equation} \label{eq:GGMSstratum}
A(\mu_\bigdot) := \bigcap_{w \in W} S_w^{\mu_w}.
\end{equation}

Fix a high weight vector $ v_{\fund_i} $ in each fundamental representation $ V_{\fund_i} $ of $ G $.  For each chamber weight $ \gamma = w \cdot \fund_i $, let $ v_\gamma = \overline{w} \cdot v_{\fund_i} $. Since $ G $ acts on $ V_{\fund_i} $,  $ G(\K) $ acts on $ V_{\fund_i} \otimes \K $.

For each $ \gamma \in \Gamma $ define the function $ \D_\gamma $ by:
\begin{equation} \label{eq:Vdef}
\begin{aligned} 
\D_\gamma: \Gr & \rightarrow \Z \\
[g] &\mapsto \val ( g \cdot v_\gamma)
\end{aligned}
\end{equation}

The functions $ \D_\gamma $ have a simple structure with respect to the semi-infinite cells.  To see this note that if $ \gamma = w \cdot \fund_i $, then $ v_\gamma $ is invariant under $ N_w $.  This immediately implies the following lemma.
\begin{Lemma}[\cite{me}] \label{th:SmDg}
Let $ w \in W $.
The function $ D_{w \cdot \fund_i} $ takes the constant value $ \langle \mu, w \cdot \fund_i \rangle $ on $ S_w^\mu $.  In fact,
\begin{equation*}
S_w^\mu = \{ L \in \Gr : D_{w \cdot \fund_i} (L) = \langle \mu, w \cdot \fund_i \rangle \textrm{ for all } i \}.
\end{equation*}
\end{Lemma}

Let $ M_\bigdot $ be a collection of integers.  Then we consider the joint level set of the functions $ D_\bigdot$,
\begin{equation}
A(M_\bigdot) := \{ L \in Gr : D_\gamma(L) = M_\gamma \text{ for all } \gamma \}.
\end{equation}

Lemma \ref{th:SmDg} shows that if $ \mu_\bigdot $ is related to $ M_\bigdot $ by (\ref{eq:mufromM}), then $ A(\mu_\bigdot) = A(M_\bigdot)$. 

Let $ \mu_1, \mu_2 $ be coweights with $ \mu_1 \le \mu_2 $.  A component of $ \overline{S_e^{\mu_1} \cap S_{w_0}^{\mu_2}} $ is called an \textbf{MV cycle} of coweight $(\mu_1, \mu_2) $.  

\begin{Theorem}[\cite{me}] \label{th:BZcycle}
Let $ M_\bigdot $ be a BZ datum of coweight $(\mu_1, \mu_2) $.  Then $ \overline{A(M_\bigdot)} $ is an MV cycle of coweight $(\mu_1, \mu_2) $.  Moreover, all MV cycles arise this way.  

Conversely, if $ Y = A(M_\bigdot) $ is an MV cycle, then 
\begin{equation*}
 \Phi(Y) := \conv \big( \{ \mu \in \cwl : t^\mu \in Y \} \big) = P(M_\bigdot).
 \end{equation*}
 
\end{Theorem}

Thus we get a coweight preserving bijection between MV cycles and polytopes, where one direction is to take an MV polytope to the closure of the associated GGMS stratum and where the other direction is to take the \textbf{moment map image} $ \Phi(Y) $ of the MV cycle $ Y $.

In general if $ Y \subset  X $ is irreducible and $ f : X \rightarrow S $ is a constructible function, then there exists a dense constructible subset $ U $ of $ Y $ such that $ f $ is constant on $ U $.  In this situation, the value of $ f $ on $ U $ is called the \textbf{generic value} of $ f $ on $ Y $.  

Using this language, Theorem \ref{th:BZcycle} says that if $ A $ is an MV cycle and $ M_\gamma $ is the generic value of $ \D_\gamma $ on A, then $ M_\bigdot $ is a BZ datum.

Note that $ \cwl $ acts on $ \Gr $ by $ \nu \cdot L := L \cdot t^\nu $.  Since $ T $ normalizes $ N_w $, we see that $ \nu \cdot S_w^\mu = S^{\mu+\nu}_w $.  So if $ A $ is a component of $ \overline{S_e^{\mu_1} \cap S_{w_0}^{\mu_2}} $, then $ \nu \cdot A $ is a component of $ \overline{S_e^{\mu_1 + \nu} \cap S_e^{\mu_2 + \nu}} $.  So $ \cwl $ acts on the set of all MV cycles. The orbit of an MV cycle of coweight $ (\mu_1, \mu_2) $ is called a \textbf{stable MV cycle} of coweight $\mu_1 - \mu_2 $.  Note that a stable MV cycle of coweight $ \mu $ has a unique representative of coweight $(\nu + \mu, \nu) $ for any coweight $\nu $.  Let $ \mathcal{M} $ denote the set of stable MV cycles and let $ \mathcal{M}(\mu) $ denote the set of those of coweight $\mu $. 

\subsection{J-relative MV cycles} \label{se:restrict}

Now we introduce the Braverman-Finkelberg-Gaitsgory crystal structure on the set of MV cycles.  Their key idea is to write each MV cycle as a pair consisting of a MV cycle for a rank 1 Levi subgroup and an MV cycle relative to the corresponding standard parabolic subgroup.  It will be convenient to work more generally (as they do) and consider such a decomposition with respect to any standard parabolic subgroup, or equivalently with respect to a subset $ J \subset I $.

Fix a subset $ J \subset I $.  Let $ \Delta_J \subset \Delta $ be the root subsystem generated by the simple roots $ \{ \alpha_j : j \in J \} $.  Let $ P^J $ be the standard parabolic subgroup of $ G $ whose Lie algebra consists of all root spaces $ \mathfrak{g}_\alpha $ with $ \alpha \in \Delta^+ \cup \Delta_J $.  Let $ N^J $ be the nilpotent radical of $ P^J $. 

Let $ G_J = P^J / N^J $ be the Levi group corresponding to $ J $.  Note that the maximal torus $ T $ of $ G $ maps isomorphically into this quotient to become a maximal torus of $ G_J $.  Also, let $ N^-_J $ be the subgroup of $ G $ whose Lie algebra contains all weight spaces $ \mathfrak{g}_\alpha $ for $  \alpha \in \Delta_J^- $.  Then $ N^-_J \subset P^J $ and maps isomorphically into $ G_J $ under the quotient.  The image of $ B $ under this quotient is denoted $ B_J $ and its unipotent radical is denoted $ N_J $.  

There is a canonical splitting $ G_J \rightarrow P^J $ of this quotient.  Under this splitting, the image of $ N_J $ is the subgroup of $ G $ (also denoted $N_J $) whose Lie algebra contains all weight spaces $ \mathfrak{g}_\alpha $ for $ \alpha \in \Delta_J^+ $.  Let $ W_J $ denote the Weyl group of $ G_J $.  Under this splitting, it is identified with the subgroup of $ W $ generated by $ \{ s_j : j \in J \} $.  Let $ \Jl \in W $ denote the long element of $ W_J $.

Let $ \Gr_J $ be the affine Grassmannian of $ G_J$.  For $ w \in W_J $ we have the unipotent subgroup $N_{w, J} = w N_J w^{-1}$.  Note that $ N_{\Jl, J} = N^-_J $.  For $ w \in W_J$ and $ \mu \in \cwl $, let $ S^\lambda_{w,J} $ denote the corresponding semi-infinite cell in $ \Gr_J$.  

For each $ \lambda \in \cwl $, we have the \textbf{restriction map} $ q_J^\lambda : S_\Jl^\lambda \rightarrow \Gr_J $ which on the level of $ \C $-points is given by $ q_J^\lambda(t^\lambda \cdot g) = t^\lambda \cdot \bar{g} $, 
where $ \bar{g} $ denotes the image of $ g \in N_v(\K) $ in the quotient $ G_J(\K) $.  These $ q_J^\lambda $ do not patch together to define a morphism from $ \Gr $ to $\Gr_J$ (see \cite{BG}).  However, as a slight abuse of notation we will usually drop by the superscript on the $ q_J $ since $\lambda $  will be obvious from the context.
 
\begin{Proposition} \label{th:imageS}
Let $ w \in W_J$.  Then
\begin{equation*}
q_J(S^\mu_w \cap S^\lambda_\Jl) = S^\mu_{w,J} \cap S^\lambda_{\Jl,J}.
\end{equation*}
\end{Proposition}

\begin{proof}
Let $ R \in S^\lambda_w $.  Then $ R = [t^\lambda g] $  for some $ g \in N_\Jl(\K) $.  Also $ R = [t^\mu h] $ for some $ h \in N_w(\K) $.  Hence there exists $ r \in G(\mathcal{O}) $ such that $ r t^\lambda g = t^\mu h$.  But we see that $ r \in P(\K) $ since all the rest of the factors in this equation are.  So we see that $ \bar{r} t^\lambda \bar{g} = t^\mu \bar{h} $ with $ r \in G_J(\mathcal{O})$.  Hence $ q_J(R) = [t^\lambda \bar{g}] = [t^\mu \bar{h}] \in S^\mu_{w, J} $ as desired.

To see that this map is onto, we use the splitting $ G_J \rightarrow G $ to provide a splitting for $ q^\lambda_J $.
\end{proof} 

For the remainder of this section $ \mu_1, \lambda, \mu_2 $ will denote coweights such that $ \mu_1 \le \lambda $, $ \lambda \le_v \mu_1 $, $\lambda \le_v \mu_2 $, $\mu_2 \le_{w_0} \lambda $.  These conditions are equivalent to $ S_e^{\mu_1} \cap S_v^\lambda \cap S_{w_0}^{\mu_2} \ne \emptyset $.  We will abbreviate these inequalities as \begin{equation*}
 \mu_1 \le_e^v \lambda \le_v^{w_0} \mu_2.
 \end{equation*}

When restricted appropriately, the map $ q_J $ is the projection of a product.
\begin{Lemma}[\cite{BFG}] \label{th:prod}
Let $ \lambda \le_v^{w_0} \mu_2 \in \cwl $.  Then, for any $ L \in S^\lambda_{\Jl, J} $, the map
\begin{align*}
S^\lambda_{\Jl, J} \times (q_J^{-1}(L) \cap S_{w_0}^{\mu_2}) &\rightarrow S^\lambda_{\Jl} \cap S_{w_0}^{\mu_2} \\
(L \cdot g, R) &\mapsto R \cdot g  \ \text{ where } g \in N^-_J(\K)  
\end{align*}
is an isomorphism.
\end{Lemma}

In other words, the action of $ N_J^-(\K) $ identifies the fibres of $ q_J : S^\lambda_{\Jl} \cap S_{w_0}^{\mu_2} \rightarrow S_{\Jl, J}^\lambda $ in a consistent manner.  So if $ L_1, L_2 \in S^\lambda_{\Jl, J}, g \in N_J^-(\K) $ with $ L_1 \cdot g = L_2$, then $ q_J^{-1}(L_1) \cap S_{w_0}^{\mu_2} $ and $ q_J^{-1}(L_2) \cap S_{w_0}^{\mu_2} $ are isomorphic by the action of $ g $ and this isomorphism does not depend on the choice of $ g $.

Note that $ q_J^{-1}(t^\lambda) \subset S^\lambda_e $ by Proposition \ref{th:imageS}.  So we see that $ q_J^{-1}(t^\lambda) \cap S^{\mu_2}_{w_0} \subset S^\lambda_e \cap S^{\mu_2}_{w_0} $ and hence is finite dimensional. Hence all the fibres are finite dimensional.

Let $ L_1, L_2 \in S^\lambda_{\Jl, J} $.  We say that two components $ B_1 $ of $q_J^{-1}(L_1) \cap S^{\mu_2}_{w_0} $ and $ B_2 $ of $ q_J^{-1}(L_2) \cap S^{\mu_2}_{w_0} $ are called \textbf{base-equivalent} if there exists $ g \in N^-_J(\K) $ such that $ L_1 \cdot g = L_2 $ and $ B_2 = B_1 \cdot g $.  A base-equivalence class of such components is called a $J$-\textbf{relative MV cycle} of coweight $ (\lambda, \mu_2) $.  

Let $ L \in S_{v, J}^\lambda $.  By the above lemma, no two components of $ q_J^{-1}(L) \cap S_{w_0}^{\mu_2} $ are base-equivalent.  So the set of components of $ q_J^{-1}(L) \cap S_{w_0}^{\mu_2} $ is in bijection with the set of $J$-relative MV cycles of coweight $(\lambda, \mu_2) $.

\begin{Proposition} \label{th:combij}
Let $ \mu_1 \le_e^v \lambda \le_v^{w_0} \mu_2 $.  The following map is a bijection:
\begin{align*}
\Comp(S^{\mu_1}_e \cap S^\lambda_{\Jl} \cap S_{w_0}^{\mu_2}) &\rightarrow \textrm{$G_J$-MV cycles of cwt $(\mu_1, \lambda) $} \times \textrm{$J$-rel. MV cycles of cwt $(\lambda, \mu_2)$}  \\
Z &\mapsto (q_J(Z), [Z \cap q_J^{-1}(L)]), \textrm{ for any } L \in q_J(Z)
\end{align*}
\end{Proposition}

\begin{proof}
By Proposition \ref{th:imageS},  $ R \in S^{\mu_1}_e \cap S^\lambda_{\Jl} $ if and only if $ q_J(R) \in S^{\mu_1}_{e, J} \cap S^\lambda_{\Jl, J} $.  Hence under the isomorphism in Lemma \ref{th:prod}, for any $ L \in  S^{\mu_1}_{e,J} \cap S^\lambda_{{\Jl}, J} $, 
\begin{equation*}
(S^{\mu_1}_{e, J} \cap S^\lambda_{\Jl, J}) \times (q_J^{-1}(L) \cap S^{\mu_2}_{w_0}) \simeq 
S^{\mu_1}_e \cap S^\lambda_{\Jl} \cap S_{w_0}^{\mu_2}.
\end{equation*}

If $ X, Y $ are varieties, then the map
\begin{equation*}
 \Comp(X) \times \Comp(Y) \rightarrow \Comp(X \times Y)  \qquad (A, B) \mapsto A \times B 
\end{equation*}
is a bijection, whose inverse is given by
\begin{equation*}
Z \mapsto (\pi_1(Z), Z \cap \pi_1^{-1}(a))
\end{equation*}
where $ \pi_1 : X \times Y \rightarrow X $ is the projection, $ a \in \pi_1(Z) $, and where we identify $ \pi_1^{-1}(a) $ with $ Y $.  

This is our situation, with $ X = S_{e, J}^{\mu_1} \cap S^\lambda_{\Jl, J}, \ X \times Y = S^{\mu_1}_e \cap S^\lambda_{\Jl} \cap S_{w_0}^{\mu_2},$ and $ \pi_1 = q_J $, except that we do not have a ``$Y$'', just a consistent way of identifying the fibres of $q_J $.
\end{proof}

Suppose that $ L \in S^\lambda_{\Jl, J}$ and $ Z $ is a component of $ q_J^{-1}(L) \cap S^{\mu_2}_{w_0}$.  Let $\nu \in \cwl$.  So $ L \cdot t^\nu \in S^{\lambda+\nu}_{\Jl, J}$ and $ A \cdot t^\nu $ is a component of $q_J^{-1}(L \cdot t^\nu) \cap S^{\mu_2+\nu}_{w_0} $.  This commutes with base-equivalence and so $\cwl$ acts on the set of $J$-relative MV cycles.  The orbit of a $J$-relative MV cycle of coweight $(\lambda, \mu_2)$ is called a \textbf{stable} $J$-\textbf{relative MV cycle} of coweight $\mu_2- \lambda $.  The set of stable $ J$-relative MV cycles of coweight $\mu $ is denoted $ \mathcal{M}^J(\mu_2 - \lambda) $.

Suppose that $ Z $ is a MV cycle of coweight $ (\mu_1, \mu_2 )$.  Then there is a unique $ \lambda $ such that $ Z \cap S_v^\lambda $ is dense in $ Z $.  In fact $ \lambda = \mu_v $ where $ \mu_\bigdot $ denotes the GGMS datum of $ Z $.  Moreover, $ Z \cap S_e^{\mu_1} \cap S_v^\lambda \cap S_{w_0}^{\mu_2} $ is a component of $ S_e^{\mu_1} \cap S_v^\lambda \cap S_{w_0}^{\mu_2} $.  Actually, the methods of \cite{me} show that any component of $ S_e^{\mu_1} \cap S_v^\lambda \cap S_{w_0}^{\mu_2} $ arises this way from a unique choice of MV cycle $ Z $.  This is because if we fix a reduced word $ \wi $ for $ w_0 $ such that $ w_k^\wi = v $, then 
\begin{equation*}
S_e^{\mu_1} \cap S_v^\lambda \cap S_{w_0}^{\mu_2} = \bigcup A^\wi(n_\bigdot)
\end{equation*}
where the union is over all $ \wi$-Lusztig datum $n_\bigdot $ of coweight $ \mu_2 - \mu_1 $ such that $ \sum_{l=1}^k n_l \beta_l^\wi = \lambda - \mu_1 $.  Theorem 5.2 in \cite{me} showed that the closure of each of these pieces is an MV cycle.  

So we have a bijection
\begin{equation*}
\text{MV cycles of coweight } (\mu_1, \mu_2) \rightarrow \bigcup_{\mu_1 \le_e^v \lambda \le_v^{w_0} \mu_2} \Comp(S_e^{\mu_1} \cap S_v^\lambda \cap S_{w_0}^{\mu_2}).
\end{equation*}

Thus we get the following.

\begin{Theorem}[\cite{BFG}] \label{th:decomp}
Taking orbits under $ \cwl$, the above argument and previous proposition gives us a bijection
\begin{align*}
\mathcal{M}(\mu) \rightarrow \bigcup_{0 \le_e^v \lambda \le_v^{w_0} \mu} \mathcal{M}_J(\lambda) \times \mathcal{M}^J(\mu-\lambda).
\end{align*}
where $ \mathcal{M}_J $ is the set of MV cycles for $ G_J$.
\end{Theorem}

We are supposed to think of this map as being the decomposition of $ \mathcal{M} $ with respect to the Levi subgroup $G_J^\vee $ of $ G^\vee$.

\subsection{Crystal structure} \label{se:BFGcrys}

We are ready to define the BFG crystal structure.  First, suppose that $ G $ is a rank 1 group with $ I = \{j\} $.  Then $\mathcal{M}(\lambda) $ consists of one element for any $ \lambda \ge 0 $.  So it is trivial to define crystal operators
\begin{equation*}
e_j : \mathcal{M}(\lambda) \rightarrow \mathcal{M}(\lambda + \alpha_j), \ f_j : \mathcal{M}(\lambda) \rightarrow \mathcal{M}(\lambda - \alpha_j) 
\end{equation*}
by setting $ e_j, f_j $ to be the obvious maps, except in the case that $ \lambda + \alpha_j \nleq 0 $, in which case we set $ e_j(A) = 0 $.

More explicitly, the only MV cycle of coweight $ (\mu_1, \mu_2) $ is $ \overline{S^{\mu_1}_e \cap S^{\mu_2}_{w_0}} $ and we define $ f_j \cdot \overline{S^{\mu_1}_e \cap S^{\mu_2}_{w_0}} := \overline{S^{\mu_1 - \alpha_j^\vee}_e \cap S^{\mu_2}_{w_0}} $.

Now that we have defined the crystal structure for rank 1 groups, we can extend it to general $ G $.  Let $ j \in I$.  Consider the decomposition of Theorem \ref{th:decomp} with $ J = \{ j \} $.  Then $ G_J $ is a rank 1 group.  So we have a well-defined crystal structure on stable MV cycles for $ G_J$.  Then, we define $ f_j, e_j $ on stable MV cycles for $ G $ by defining them to act on the first factor in the decomposition of Theorem \ref{th:decomp}. 

Braverman-Finkelberg-Gaitsgory proved that this defines a crystal structure which is isomorphic to the crystal for the Verma module and hence it is isomorphic to the LBZ crystal structure on the set of MV polytopes.  The goal of the rest of this section is to prove that the map $ Y \mapsto \Phi(Y) $ achieves this isomorphism.
\begin{Theorem} \label{th:LBZ=BFG}
Under the bijection between MV polytopes and MV cycles, the LBZ and BFG crystal structures agree.
\end{Theorem}

\subsection{Relative chamber weights}

It is necessary to understand better the $ J$-relative MV cycles.  We say that $ \gamma $ is a $J$-\textbf{relative chamber weight} if $ \gamma = w \cdot \fund_i $ for some $ i \in I  $ and for some $ w $ which is a maximal length representative of the left coset $ W_J w $.  Let $ \Gamma^J $ denote the set of $J$-relative chamber weights.

\begin{Lemma} \label{th:relequiv}
Let $ \gamma $ be a $J $-relative chamber weight.  Then the following hold.
\begin{enumerate}
\item $ \langle \alpha_j^\vee, \gamma \rangle \le 0 $ for all $ j \in J $.
\item If $ g \in N^-_J$, then $ g v_\gamma = v_\gamma $.
\end{enumerate}
\end{Lemma}

\begin{proof}
(i) Since $ \gamma \in \Gamma^J$, $ \gamma = w \cdot \fund_i $ for some $ i \in I $ and some $ w $ which is a maximal length representative of the left coset $ W_J w $.  So if $ j \in J$, then $ s_j w < w $.  Hence $ w^{-1} \cdot \alpha_j^\vee $ is a negative coroot.  So
\begin{equation*}
\langle w^{-1} \cdot \alpha_j^\vee, \fund_i \rangle \le 0  \Rightarrow \langle \alpha_j^\vee, w \cdot \fund_i \rangle \le 0 
\end{equation*}
as desired.

(ii)  It suffices to show that $ E_j \cdot v_\gamma = 0 $ for all $ j \in J$, since these generate the Lie algebra of $ N^-_J $ (here $ E_j $ denotes a non-zero element of the $-\alpha_j $ weight space of $ \mathfrak{g} $).  

Consider the $ j$th copy of $ \mathfrak{sl}_2 $ inside $ \mathfrak{g} $ acting on $ V_{\fund_i}$.  Since $ \gamma $ is a Weyl orbit of the highest weight, $ v_\gamma $ is a extremal vector for this action.  However, $ \langle \alpha_j^\vee, \gamma \rangle \le 0 $, so it cannot be a highest weight vector.  Hence it is a lowest weight vector and so $ E_j \cdot v_\gamma = 0 $.
\end{proof}

\begin{Proposition} \label{th:relval}
 Let $ \gamma $ be a $J$-relative chamber weight.  Suppose that $ B_1, B_2 $ are two base-equivalent $J$-relative MV cycles.   Then the generic value of $ \D_\gamma $ on $ B_1 $ equals the generic value of $ \D_\gamma $ on $ B_2 $.

Moreover, let $ Z $ be a component of $ S^{\mu_1}_1 \cap S^\lambda_{\Jl} \cap S^{\mu_2}_{w_0} $ and let $ L \in q_J(Z) $.  Then the generic value of $ \D_\gamma $ on $ Z $ equals the generic value of $ \D_\gamma $ on $ Z \cap q_J^{-1}(L) $.
\end{Proposition}

\begin{proof}
Since $ B_1, B_2 $ are base-equivalent, there exists $ g \in N^-_J(\K) $ such that $ B_1 \cdot g = B_2 $.  So it suffices to show that $ \D_\gamma(R) = \D_\gamma( R \cdot g) $ for all $ R \in \Gr $.

Suppose that $R = [h] $ and recall that $ \D_\gamma([h]) = \val (h \cdot v_\gamma) $.  So $ \D_\gamma ([h] \cdot g) = \val( hg \cdot v_\gamma) $.  By the second part of Lemma \ref{th:relequiv}, we see that $ g \cdot v_\gamma = v_\gamma $ and hence the result follows.

For the second part, note that if $ L_1, L_2 \in q_J(Z) $, then $ q_J^{-1}(L_1), q_J^{-1}(L_2) $ are base equivalent.  Hence the first half of this proposition shows that the generic value of $\D_\gamma $ on $ Z \cap q_J^{-1}(L) $ does not depend on $ L \in q_J(Z) $.  Let $ k $ denote this value.  So for each $ L \in q_J(Z)$, we have a dense set $U_L $ of $ Z \cap q_J^{-1}(L) $ on which $ \D_\gamma $ equals $ k$.  Since each $ U_L $ is dense in the fibre over $ L $,  $ U = \cup_L U_L $ is dense in $ Z $.  So we have an dense set $ U $ in $ Z $ on which $ \D_\gamma $ assumes the value $ k $.  Hence the generic value of $ \D_\gamma $ on $ Z $ is $ k $.
\end{proof}

We are now in a position to complete our proof.

\begin{proof}[Proof of Theorem \ref{th:LBZ=BFG}]
Let $ \mu_1 \le \mu_2$ and let $ Y$ be an MV cycle of coweight $ (\mu_1, \mu_2) $.  So by Theorem \ref{th:BZcycle}, $ Y = \overline{A(M_\bigdot)}$ for some BZ datum $ M_\bigdot $.  Let $ \mu_\bigdot $ be the associated GGMS datum.  

Let $ Y' = f_j \cdot Y $ be the image of $ Y $ under the BFG crystal operator.  Let $ M'_\bigdot $ and $ \mu'_\bigdot $ denote its BZ and GGMS data.  We would like to show that $ f_j \cdot P(M_\bigdot) = P(M'_\bigdot) $ where $ f_j $ denotes the LBZ crystal operator.  By Theorem \ref{th:LBZ}, it suffices to prove that if $ \gamma $ is a $j$-relative chamber weight, then $ M'_\gamma = M_\gamma $.

We will use the above results with $ J = \{j \} $ and will abbreviate $ \{ j \} $ by $ j $.  Note that $ \Jl = s_j $.  Let $ \lambda = \mu_{s_j}$.  Since $ A(M_\bigdot) = \cap_w S^{\mu_w}_w $, we see that $ Z := Y \cap S^{\mu_1}_e \cap S^\lambda_{s_j} \cap S^{\mu_2}_{w_0}$ is dense in $Y $ .  

Let $ A = q_j(Z) $.  Choose $ L \in A $ and let $ B = Z \cap q_j^{-1}(L) $.  So $ Z $ corresponds to $ (A, [B]) $ in Proposition \ref{th:combij}.  We know that $ A $ is a component of $  S^{\mu_1}_{1,j} \cap S^\lambda_{s_j, j} $ which is irreducible, so $ A = S^{\mu_1}_{1,j} \cap S^\lambda_{s_j, j}$.

Let $ A' = f_j \cdot A = S^{\mu_1-\alpha_j^\vee}_{1,j} \cap S^\lambda_{s_j,j}$.  Fix $ L' \in S^{\mu_1 - \alpha_j^\vee}_{1,j} \cap S^\lambda_{s_j,j} $.  So $ L' = L \cdot g $ for some $ g \in N^-_j(\K) $.  Then the definition of the BFG crystal structure says that $ Z' := Y' \cap S^{\mu_1- \alpha_j^\vee}_e \cap S^\lambda_{s_j} \cap S^{\mu_2}_{w_0}  $ corresponds to the pair $(A', [B \cdot g]) = (A', [B]) $ in the bijection of Proposition \ref{th:combij}. 

Let $ \gamma $ be a $j$-relative chamber weight.  The generic value of $ \D_\gamma $ on $ Z'$ equals the generic value of $ \D_\gamma $ on $ B \cdot g $ by Proposition \ref{th:relval}.  Also, $ B \cdot g$ and $ B $ are base-equivalent, so also by Proposition \ref{th:relval}, the generic value of $ \D_\gamma $ on $ B \cdot g $ equals the generic value of $ \D_\gamma $ on $ B$.  Finally another application of Proposition \ref{th:relval} shows that the generic value of $ \D_\gamma $ on $ B $ equals the generic value of $ \D_\gamma $ on $ Z $.  So the generic value of $ D_\gamma $ on $ Z' $ equals the generic value of $ D_\gamma $ on $Z$. Since $ Z', Z $ are dense in $ Y', Y $, we conclude that $ M'_\gamma = M_\gamma $ as desired.
\end{proof}

\section{Anderson-Mirkovi\'c conjecture}  \label{se:AM}

\subsection{The conjecture}

Anderson-Mirkovi\'c gave a conjectural description of a crystal structure on MV polytopes.  We learnt of this conjecture from personal communication with Anderson in May 2003.  We begin by stating their conjecture.

Fix $ j \in I $ and let $ W^- := \{ w \in W : s_j w < w \}, W^+ := \{ w \in W : s_j w > w \}$.  

Let $ P = \conv(\mu_\bigdot) $ be an MV polytope.  Let $ c = M_{\fund_j} - M_{s_j \cdot \fund_j} - 1 $ and let $ r : \realt \rightarrow \realt $  be defined by $ r(\alpha) = s_j \cdot \alpha + c \alpha_j^\vee $.  This reflection takes $ \mu_{s_j} $ to $ \mu_e - \alpha_j^\vee $.   Now we define
$ P' = AM_j \cdot P = \conv(\mu'_\bigdot) $ to be the smallest pseudo-Weyl polytope such that
\begin{enumerate}
\item $ \mu'_w = \mu_w $ for all $ w \in W^-$,
\item $ \mu'_e = \mu_e - \alpha^\vee_j $,
\item $P'$ contains $ \mu_w $ for all $ w \in W^+ $,  and
\item if $w \in W^- $ is such that $ \langle \mu_w, \alpha_j \rangle \ge c $, then $ P' $ contains $ r ( \mu_w) $. 
\end{enumerate}
Note that if $ w \in W^- $ and $ \langle \mu_w, \alpha_j \rangle < c$, then $ r(\mu_w) \ngeq_w \mu_w $ and so there is no way for $ r(\mu_w) $ to be in $ P'$. 

Hence the vertices that stay the same in this crystal structure are exactly those that stay the same in the LBZ=BFG crystal structure.  The rest of the vertices are adjusted so as to contain both the old polytope and the images of some of the fixed vertices under $ r $.  

\begin{Conjecture}[Anderson-Mirkovi\'c] \label{th:AM}
For any MV polytope $ P $ and any $ j \in I $, $ AM_j \cdot P $ is an MV polytope.  Moreover $ AM_j \cdot P = f_j \cdot P $ where $ f_j $ denotes the BFG crystal operator.
\end{Conjecture}

The definition of $ AM_j $ only asks for the smallest pseudo-Weyl polytope and so does not make any reference to the tropical Pl\"ucker relations.  Thus, in the cases where the conjecture holds, we get a more explicit picture than for the LBZ crystal structure.  Moreover, the conjecture gives a simple method for inductively generating all of the MV polytopes without reference to the tropical Pl\"ucker relations.  This was one of the original motivations of Anderson-Mirkovi\'c.

We will show that the conjecture holds for a class of Lie algebras including $ \mathfrak{sl}_n$, but that it fails for $ \mathfrak{sp}_6 $.

\begin{Example} \label{eg:AM}
We continue from Example \ref{eg:LBZ}.  In this case we see that $ c = -1 $, so the reflected vertices are $ r(0,-1,1) = (-2,1,0),\ r(1,-1,0) = (-2,2,0),\ r(1,0,-1) = (-1,2,-1) $.  The minimal pseudo-Weyl polytope containing these vertices is the polytope we previously calculated. \samepage
\begin{center}
\epsfig{file=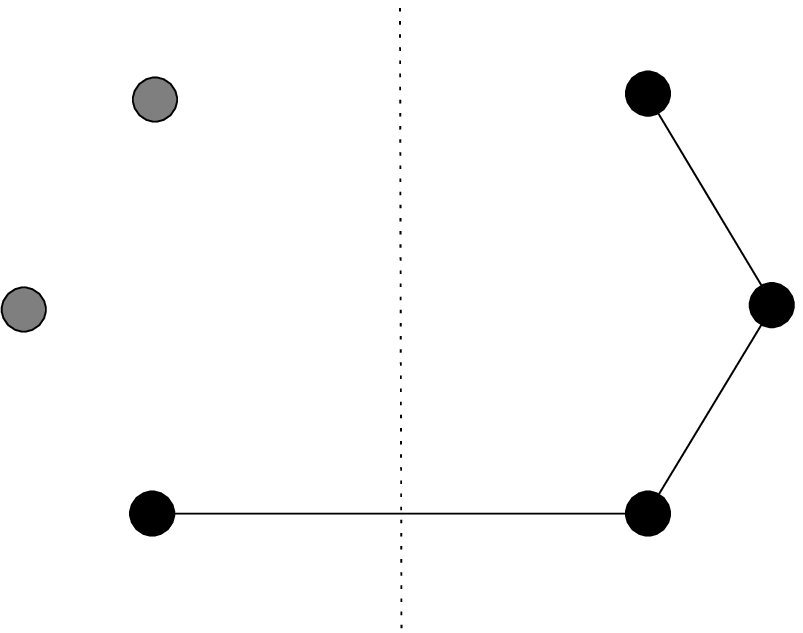,height=2cm}
\end{center}
\end{Example}

\subsection{BZ data and the AM conjecture}

We begin by rephrasing the AM conjecture to fit with the language of BZ datum.
Let $ M_\bigdot $ be a BZ datum and let $ j \in I$.  Let $ \Gamma^j = \cup_i W^- \cdot \fund_i $ be the set of $j$-relative chamber weights and let $ \Gamma_j = \Gamma \setminus \Gamma^j $. Let $ M'_\bigdot $ be defined by
\begin{equation} \label{eq:M'}
M'_\gamma = 
\begin{cases} M_\gamma, \ \text{ if } \gamma \in \Gamma^j \\
 \min\big(M_\gamma, M_{s_j \cdot \gamma} + c \langle  \alpha_j^\vee, \gamma \rangle \big), \  \text{ if } \gamma \in \Gamma_j.
\end{cases}
\end{equation}
Here $ c $ is as in the definition of $ AM_j $.

Note that by this definition $ M'_{\fund_j} = \min(M_{\fund_j},M_{\fund_j} - 1) = M_{\fund_j} - 1$ .

\begin{Proposition} \label{th:AMBZ}
Suppose that $ M'_\bigdot $ satisfies the edge inequalities.  Then $ P(M'_\bigdot) = AM_j \cdot P(M_\bigdot) $.
\end{Proposition} 

\begin{proof}
Since $ M'_\bigdot $ satisfies the edge inequalities, $P' := P(M'_\bigdot) $ is a pseudo-Weyl polytope.  So to prove this proposition, we must show that $P' $ satisfies condition (i)-(iv) in the definition of $ AM_j$ and that it is the minimal pseudo-Weyl polytope satisfying those conditions.  

It is clear that $P'$ satisfies conditions (i) and (ii), so now we consider (iii), (iv).  Let $\mu'_\bigdot $ denote the GGMS datum of $ P' $.  In general to show that $ \nu \in P' $, we will show that $ \nu \ge_v \mu_v' $ for all $ v \in W$.  This is how we will establish (iii) and (iv) (with the role of $ \nu $  played by $ \mu_w, r(\mu_w) $ respectively).

Let $ w \in W^+$ and let $ v \in W$.  Then for any $ \gamma = v \cdot \fund_i $, $M'_\gamma \le M_\gamma$.  Hence $ \mu_v' \le_v \mu_v $.  Since $ P $ is a pseudo-Weyl polytope, $ \mu_w \ge_v \mu_v $.  So $ \mu_w \ge_v \mu'_v $ for all $ v \in W $.  Hence $ \mu_w \in P' $ and so condition (iii) holds.

Let $ w \in W^- $ such that $ \langle \mu_w, \alpha_j \rangle \ge c $.  Let $ d = \langle \mu_w, \alpha_j \rangle - c $.  So $ d \ge 0 $ and $ r(\mu_w) = \mu_w - d \alpha_j^\vee$.  We would like to show that $ r(\mu_w) \in P'$.  To do this we will first show that $ r(\mu_w) \ge_{s_j w} \mu_{s_j w}$, then that $ r(\mu_w) \ge_v \mu_v $ for all $ v \in W^+ $, and finally that $ r(\mu_w) \ge_v \mu_v $ for all $ v \in W^- $.

Let $ i \in I$.  Let $ \gamma = s_j w \cdot \fund_i $.  We want to show that $ \langle r(\mu_w), \gamma \rangle \ge M'_\gamma$.  We deal with two cases.

Suppose that $ \gamma \in \Gamma^j $.   Then $ \langle \alpha_j^\vee, \gamma \rangle \le 0 $ by Lemma \ref{th:relequiv}.  So 
\begin{equation*}
 \langle r(\mu_w), \gamma \rangle = \langle \mu_w, \gamma \rangle - d \langle \alpha_j^\vee, \gamma \rangle \ge \langle \mu_w, \gamma \rangle.
 \end{equation*}
But $ \langle \mu_w, \gamma \rangle \ge M_\gamma $ since $ \mu_w \ge_{s_j w} \mu_{s_j w}$.  Also $ M_\gamma = M'_\gamma$.  Hence $ \langle r(\mu_w), \gamma \rangle \ge M'_\gamma $ as desired.

Suppose that $ \gamma \in \Gamma_j $.  Then
\begin{equation*}
\langle r(\mu_w), \gamma \rangle = \langle s_j \cdot \mu_w + c \alpha_j^\vee, s_j w \cdot \fund_i \rangle = M_{w \cdot \fund_i} + c \langle \alpha_j^\vee, \gamma \rangle \ge M'_\gamma,
\end{equation*}
as desired.

So $ \langle r(\mu_w), s_j w \cdot \fund_i \rangle \ge M'_{s_j w \cdot \fund_i} $ for all $ i $.  Hence $ r(\mu_w) \ge_{s_j w} \mu'_{s_j w} $.

Now for all $ v \in W^- $, $ \mu_w \ge_v \mu_v $.  Note that $ \alpha \ge_v \beta \Leftrightarrow r(\alpha) \ge_{s_j v} r(\beta) $ and so $ r(\mu_w) \ge_{s_j v} r(\mu_v) $.  So by above $ r(\mu_w) \ge_{s_j v} \mu'_{s_j v} $  for all $ v \in W^- $.  

Also $  r(\mu_w) - \mu_w = -d \alpha_j^\vee $.  Let $ v \in W^- $.  By Lemma \ref{th:relequiv}, $\alpha_j^\vee \le_v 0 $.  So we see that $ r(\mu_w) \ge_v \mu_v = \mu'_v $ for all $ v \in W^-$.

Since $ W^- \cup s_j W^- = W $ we see that $ r(\mu_w) \ge \mu'_v $ for all $ v \in W $.  This shows that condition (iv) is satisfied.

So $ P' $ satisfies (i)-(iv). 

To see that it is minimal let $ P(M''_\bigdot) $ be  a pseudo-Weyl polytope satisfying (i)-(iv).  Then $ M''_{\fund_j} = M_{\fund_j} - 1 $ and $ M''_\gamma = M_\gamma $ for all $ \gamma \in \Gamma^j$.

By (iii), $ \mu_w \in P(M''_\bigdot)$ for all $ w \in W^+$.  So we see that $ M''_\gamma \le M_\gamma $ for all $ \gamma \in \Gamma_j $.

We would also like to show that $ M''_\gamma \le M_{s_j \cdot \gamma } + c \langle \alpha_j^\vee, \gamma \rangle $ if $ \gamma \in \Gamma_j $.  

If there exists $ w \in W^- $ such that $ w \cdot \fund_i = s_j \cdot \gamma $ and $ \langle \alpha_j^\vee, \mu_w \rangle \ge c $ then $ r(\mu_w) \in P(M''_\bigdot) $ by (iv).  So $ \langle r(\mu_w), \gamma \rangle \ge M''_\gamma $.  But 
\begin{equation*}
\langle r(\mu_w), \gamma \rangle = M_{s_j \cdot \gamma} + c \langle \alpha_j^\vee, \gamma \rangle
\end{equation*}
and so we get the desired conclusion.

If such a $ w $ does not exist, then we can find $ w \in W^- $ such that $ w \cdot \fund_i = s_j \cdot \gamma $ and $ \langle \alpha_j^\vee, \mu_w \rangle < c $.  In this case, 
\begin{equation*}
M_\gamma \le \langle \mu_w, \gamma \rangle < \langle r(\mu_w), \gamma \rangle = M_{s_j \cdot \gamma} + c \langle \alpha_j^\vee, \gamma \rangle
\end{equation*}
and so $ M''_\gamma \le M_{s_j \cdot \gamma} + c \langle \alpha_j^\vee, \gamma \rangle $ as desired.

Hence we see that $M'_\gamma \ge M''_\gamma $ for all $ \gamma$.  So $ P' \subset P(M''_\bigdot) $.  Hence $ P' $ is minimal and so $ AM_j \cdot P = P' $ as desired.
\end{proof}

\subsection{$j$-close Lie algebras}
Let $ j \in I $.  We say that $ \mathfrak{g}^\vee $ is $j$-\textbf{close} if there exists a function $ H : \Gamma_j \rightarrow \Z $, such for any chamber weight $ \gamma \in \Gamma_j$ with $\gamma \ne \fund_j $, there exists $ v \in W $ and $ i, k \in I $ such that 
\begin{gather*}
vs_k > v, \ vs_i > v, \ s_j v = vs_k, \ a_{ik} = a_{ki} = -1 , \ vs_i \cdot \fund_i = \gamma,\\ 
\text{ and } H(\gamma) > H(\delta) \text{ where } \delta = v \cdot \fund_k .
\end{gather*}

Figure \ref{fig:jclose} contains a picture of these chamber weights and Weyl group elements.  We label an edge of the hexagon by the unique chamber weight whose corresponding hyperplane contains that edge but no other vertex of the hexagon.

If the $j$-close condition holds, then the action of $ f_j $ is easy to calculate by an inductive procedure.  This is exploited in the proof of Theorem \ref{th:jclosetrue}.  We also have a definition of $j$-close for non-simply laced groups, but it is more complicated so we omit it.  

\begin{Lemma}
For any $ j $, $ \mathfrak{sl}_n $ is $j$-close.
\end{Lemma}
\begin{proof}
In this case, $ I = \{ 1, \dots, n-1 \}$ and $ W = S_n$.  Also $ \fund_i = (1, \dots, 1 ,0, \dots, 0) $ and so we identify $ W \cdot \fund_i $  with the set of  $ i $ element subsets of $\{1, \dots, n \} $.

Let $ j \in \{1, \dots, n-1 \} $.  Then $ \Gamma_j $ consists of all proper subsets $ \gamma $ of $ \{1, \dots, n \} $ such that $ j \in \gamma $ and $ j+1 \notin \gamma $ (this can be seen from Lemma \ref{th:relequiv}).

Define a function $ H $ on $ \Gamma_j $ by 
\begin{equation*}
H(\gamma) := \# \{ k \in \gamma : k > j+1 \} - \# \{ k \in \gamma : k < j \}.
\end{equation*}

Let $ \gamma \in \Gamma_j $, with $ \gamma \ne \fund_j $, and $ \# \gamma $ = i.  So $ j \in \gamma $ and $ j+1 \notin \gamma $ and $ \gamma \ne \{1, \dots, j \} $.  Hence either we can find $ m \notin \gamma $ with $ m < j $ or we can find $ m \in \gamma $ with $ m > j+1 $.

Suppose the first case holds.  Then let $ v \in S_n $ be such that $ v(\{1, \dots, i-1 \}) = \gamma \setminus \{j \}$, $ v(i) = m$, $v(i+1) = j$, and $v(i+2) = j+1 $.  Then letting $ k = i+1 $ gives a pair satisfying the condition in the definition of $ j$-close (in particular, $ \delta = v \cdot \fund_k = \gamma \cup \{ m \} $ and so $ H(\delta) = H(\gamma) - 1 $).

Suppose the second case holds.  Then let $ v \in S_n $ be such that $ v (\{1, \dots, i-2 \}) = \gamma \setminus \{j,m \}$, $v(i-1) = j $, $ v(i) = j+1$, and $v(i+1) = m$.  Then letting $k = i-1 $ gives a pair satisfying the condition in the definition of $ j$-close (in particular, $ \delta = v \cdot \fund_k = \gamma \setminus \{m \} $ and so $ H(\delta) = H(\gamma) - 1 $). 

Hence the result follows.
\end{proof}

Actually, it is possible to show that if $ \mathfrak{g}^\vee $ is simple, simply-laced, and $j$-close for some $ j$, then $ \mathfrak{g}^\vee = \mathfrak{sl}_n $ for some $ n$.

\begin{Theorem} \label{th:jclosetrue}
Suppose that $ \mathfrak{g}^\vee $ is $j$-close.  Then the Anderson-Mirkovi\'c conjecture holds for this $ j $.
\end{Theorem}

\begin{proof}
Let $ P = P(M_\bigdot) $ be an MV polytope with BZ datum $ M_\bigdot $ and GGMS datum $ \mu_\bigdot $.

Let $ M'_\bigdot $ be the BZ datum of $ P' = e_j \cdot P $.  Let $ \mu'_\bigdot $ be its GGMS datum.

By Proposition \ref{th:AMBZ}, it suffices to show that $ M'_\bigdot $ is given by equations (\ref{eq:M'}).  By Theorem \ref{th:LBZ} and the remarks after its proof, we see that $ M'_\gamma = M_\gamma $ if $ \gamma \in \Gamma^j$.  So it remains to check the case $ \gamma \in \Gamma_j $.

 We proceed by induction on the set $ \Gamma_j $ using the function $ H $.  Since $ G $ is $j$-close, if $ \gamma \in \Gamma_j $ and $ \gamma \ne \fund_j$, there exists $ \delta \in \Gamma_j $ with $ H(\delta) < H(\gamma) $.  So we see that $ H(\fund_j) <  H(\gamma) $ for all $ \gamma \in \Gamma_j $.  Hence the base case of our induction is $ \fund_j $. Since $ \mu'_e = \mu_e - \alpha_j^\vee $,  $ M'_{\fund_j} = \langle \mu'_e , \fund_j \rangle = M_{\fund_j} - 1 $ as desired.

Now, suppose that $ \gamma \in \Gamma_j $ and $ M'_\delta = \min(M_\delta, M_{s_j \cdot \delta} + c \langle \alpha_j^\vee , \delta \rangle) $ for all $ \delta \in \Gamma_j $ with $ H(\delta) < H(\gamma) $.  Since $ G $ is $j$-close there exist $v \in W $ and $ k \in I $ such that $ vs_k > v, \ vs_i > v, \ s_j v = vs_k, \ a_{ik} = a_{ki} = -1 , v s_i \cdot \fund_i = \gamma $ and $ H(\gamma) > H(\delta) $, where $ \delta = v \cdot \fund_k $.  

\begin{figure} 
\begin{center}
\psfrag{v}{$v$} \psfrag{gamma}{$\gamma$} \psfrag{delta}{$\delta$} \psfrag{vs_i}{$vs_i$} \psfrag{vs_k}{$vs_k$}
\psfrag{vs_is_k}{$vs_is_k$} \psfrag{vs_ks_i}{$vs_ks_i$} \psfrag{vs_is_ks_i}{$vs_is_ks_i$} \psfrag{v fund_i}{$v \cdot \fund_i$} \psfrag{vs_is_k fund_k}{$vs_is_k \cdot \fund_k$} \psfrag{vs_k fund_k}{$vs_k \cdot \fund_k$} \psfrag{vs_ks_i fund_i}{$vs_ks_i \cdot \fund_i$} \psfrag{sj}{$s_j$}
\epsfig{file=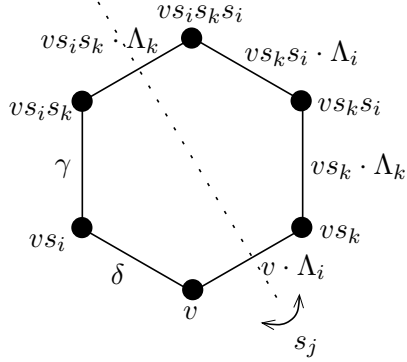,height=5cm}
\caption{The weights occurring in the definition of $ j$-close.}
\label{fig:jclose}
\end{center}
\end{figure}

Since $ M'_\bigdot $ satisfies the tropical Pl\"ucker relation at $ (v, i, k) $ (see (\ref{eq:A2trop})), we see that
\begin{equation} \label{eq:M'tpr}
M'_{vs_i \cdot \fund_i} + M'_{vs_k \cdot \fund_k} = \min( M'_{v \cdot \fund_i} + M'_{vs_i s_k \cdot \fund_k}, M'_{v s_k s_i \cdot \fund_i} + M'_{v \cdot \fund_k} )
\end{equation}

Now, $ s_j v s_k = v < v s_k $, so $ v s_k \in W^- $.  Similarly, $ v s_k s_i, v s_i s_k s_i \in W^- $ and so $ v \cdot \fund_i, v s_k \cdot \fund_k, v s_k s_i \cdot \fund_i, v s_i s_k \cdot \fund_k \in \Gamma^j $.  As already noted $ \gamma, \delta \in \Gamma_j $.  Note also that $ s_j \cdot \delta = s_j v \cdot \fund_k = v s_k \cdot \fund_k $.  Similarly, $ s_j \cdot \gamma = v s_k s_i \cdot \fund_i $.

By hypothesis, $ H(\delta) < H(\gamma) $, so by induction, $ M'_\delta = \min(M_\delta, M_{s_j \cdot \delta} + c \langle \alpha_j^\vee, \delta \rangle) $.

So (\ref{eq:M'tpr}) becomes
\begin{align*}
M'_\gamma + M_{vs_k \cdot \fund_k} &= 
\min \big( M_{v \cdot \fund_i} + M_{vs_i s_k \cdot \fund_k}, M_{v s_k s_i \cdot \fund_i} + M'_\delta \big) \\
\Rightarrow M'_\gamma + M_{vs_k \cdot \fund_k} &= 
\min \big( M_{v \cdot \fund_i} + M_{v s_i s_k \cdot \fund_k}, M_{v s_k s_i \cdot \fund_i} + 
\min(M_\delta, M_{s_j \cdot \delta} + c \langle \alpha_j^\vee, \delta \rangle) \big) \\
\Rightarrow M'_\gamma + M_{vs_k \cdot \fund_k} &= 
\min \big( M_{v \cdot \fund_i} + M_{v s_i s_k \cdot \fund_k}, M_{v s_k s_i \cdot \fund_i} + M_\delta, 
M_{s_j \cdot \gamma} + M_{s_j \cdot \delta} + c\langle \alpha_j^\vee, \delta \rangle \big) 
  \\
\Rightarrow M'_\gamma &= \min(M_\gamma, M_{s_j \cdot \gamma} + c \langle \alpha_j^\vee, \delta \rangle)
\end{align*}
where in the last step we have used that $M_\bigdot $ satisfies the tropical Pl\"ucker relation.

So it remains to show that $ \langle \alpha_j^\vee, \gamma \rangle = \langle \alpha_j^\vee, \delta \rangle $.  

Now $s_j \cdot \delta + \gamma = v s_k \cdot \fund_k + v s_i \cdot \fund_i = v \cdot \fund_k + v s_k s_i \cdot \fund_i = \delta + s_j \cdot \gamma $, where we have used $ s_k \cdot \fund_k + s_i \cdot \fund_i = \fund_k + s_k s_i \cdot \fund_i$.  Expanding out the definition of $ s_j \cdot \delta, s_j \cdot \gamma $, shows that $ \langle \alpha_j^\vee, \gamma \rangle = \langle \alpha_j^\vee, \delta \rangle$, which completes the induction step and hence the proof.
\end{proof}

As a corollary, we can now give an explicit description of the crystal structure on MV polytopes for $ SL_n $.  Recall that as usual we identify the set of chamber weights for $ \mathfrak{sl}_n $ with the set of proper subsets of $ \{1, \dots, n\} $.

\begin{Corollary} \label{th:AMtypeA}
Let $ P $ be an MV polytope for $ \mathfrak{sl}_n $.  Let $ M_\bigdot, M'_\bigdot $ denote the BZ data of $ P$, $f_i \cdot P $.  Then
\begin{equation*}
M'_\gamma = \begin{cases}
\min(M_\gamma, M_{ \gamma \setminus j \cup j+1}  + c) &\text{ if $ j \in \gamma $ and $ j+1 \notin \gamma$}\\
M_\gamma &\text{otherwise},
\end{cases}
\end{equation*}
where $ c = M_{\{1, \dots, j\}} - M_{\{1, \dots, j-1, j+1\}} -1 $.
\end{Corollary}

\subsection{Counterexample for $\mathfrak{sp}_6$} \label{se:counter}
To give a counterexample, we will exhibit a MV polytope $ P $ such that $ AM_j \cdot P $ is not an MV polytope.  More specifically, the hypothesis of Proposition \ref{th:AMBZ} will hold, so that $ AM_j \cdot P = M'_\bigdot $.  However, $ M'_\bigdot $ will not satisfy the tropical Pl\"ucker relations.

We work with $ \Lieg^\vee = \mathfrak{sp}_6 $ so that $ \realt = \R^3, \cwl = \Z^3 $ and the root datum is
\begin{gather*}
\alpha^\vee_1 = (1,-1,0), \ \alpha^\vee_2 = (0,1,-1), \ \alpha^\vee_3 = (0,0,1), \\
\alpha_1 = (1,-1,0), \ \alpha_2 = (0,1,-1), \ \alpha_3 = (0,0,2). 
\end{gather*}
The fundamental weights are
\begin{equation*}
\fund_1 = (1,0,0), \fund_2 = (1,1,0), \fund_3 = (1,1,1), 
\end{equation*}
which we abbreviate as $ 1, 12, 123 $.

The Weyl group is the set of signed permutations of $ \{1,2,3\} $ and so the chamber weights are the signed subsets of $ \{1,2,3\} $ (for example we write $ 1 - \!2 $ for the chamber weight $(1,-1,0) $). 

Fix some integer $ x \ge 2 $ and consider the MV polytope $ P $ with vertices $ (0,0,0)$, $(0,2,0)$, $(0,0,x)$, and $ (0,2,x) $.  The corresponding collection $ M_\bigdot $ is
\begin{gather*}
M_1 = M_{12} = M_{13} = M_{123} = M_{2}= M_{23} = M_3 = M_{-123} = M_{-12} = M_{-13} = M_{-1} = 0, \\
M_{1-3} = M_{12-3} = M_{2-3} = M_{-1 2 -3} = M_{-1-3} = M_{-3} = -x, \\
M_{1-2} = M_{1-23} = M_{-23} = M_{-1-23} = M_{-1-2} = M_{-2} = -2, \\
M_{1-2-3} = M_{-2-3} = M_{-1-2-3} = -x-2.
\end{gather*}

Now we consider $ j = 1 $ and define $M'_\bigdot $ by (\ref{eq:M'}).  We see that 
\begin{gather*}
M'_1 = M'_{13} = -1, \ M'_{1-23} = M'_{-23} = M'_{-2} = M'_{1-2} = -2, \\
M'_{1-3} = -x-1, \ M'_{1-2-3} = M'_{-2-3} = -x-2,
\end{gather*}
where we have just listed those $ \gamma \in \Gamma_1$.

These $ M'_\bigdot $ satisfy the edge inequalities, so by Proposition \ref{th:AMBZ}, $AM_1 \cdot P = P(M'_\bigdot) $.  We note that $ P(M'_\bigdot) $ has vertices $ (0,0,0)$, $(-1,1,0)$, $(0,2,0)$, $(0,0,x)$, $(0,2,x)$, $(-1,1,x) $.

\begin{figure}
\begin{center}
\psfrag{sj}{$s_j$}
\epsfig{file=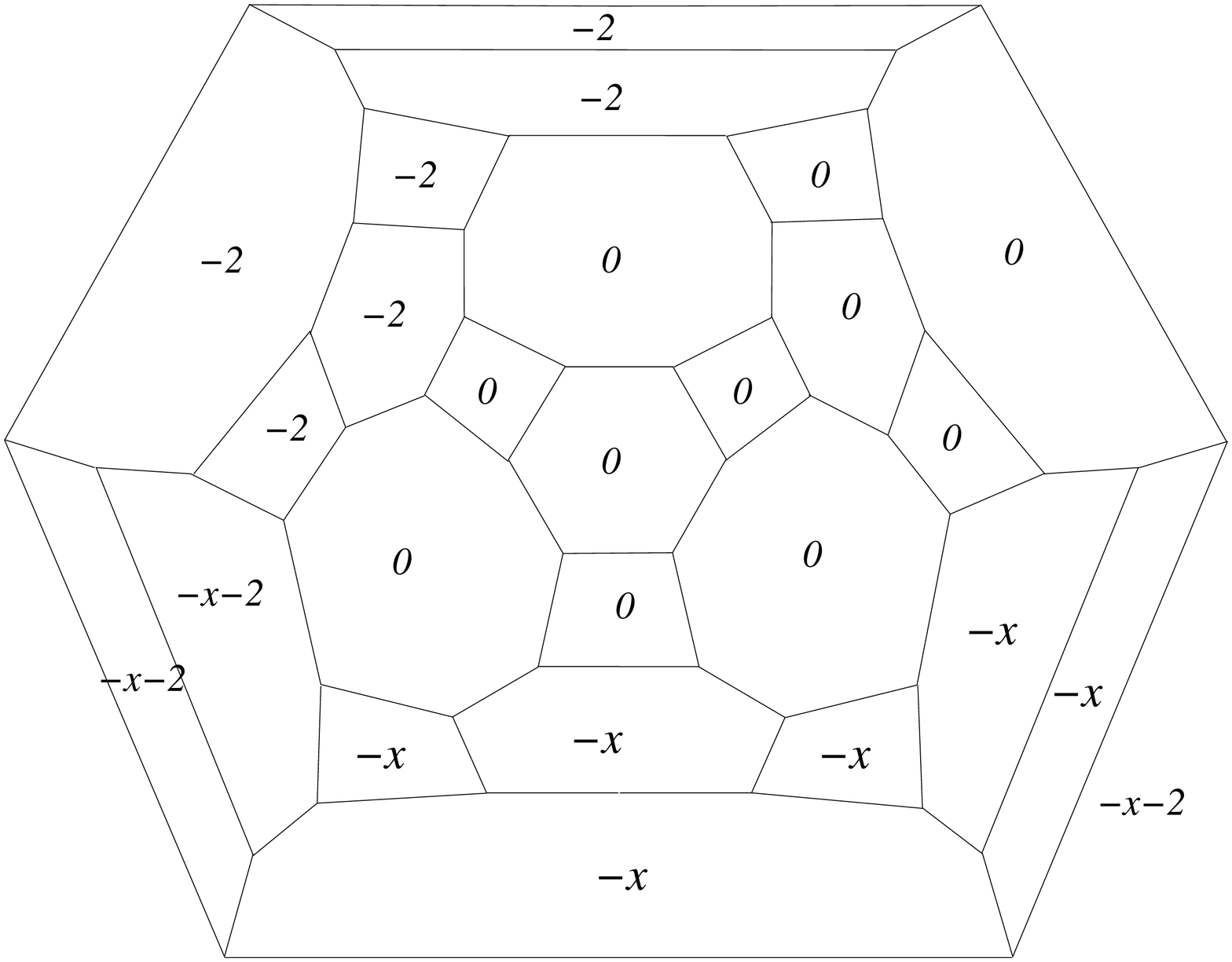,height=6cm} 
\epsfig{file=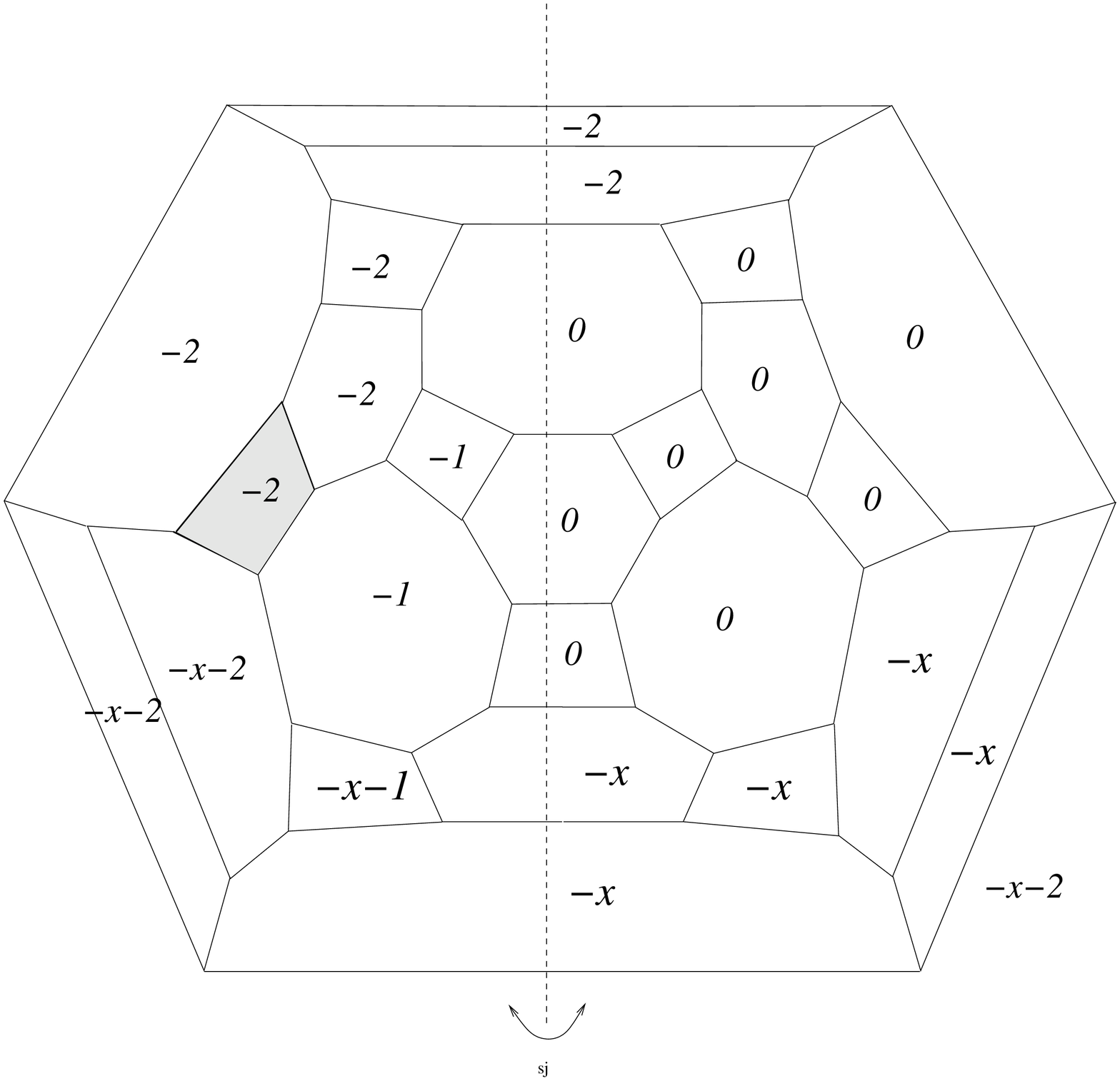,height=7cm}
\caption{The polytopes $ P $ and $ AM_j \cdot P $.}
\label{fig:M}
\end{center}
\end{figure}

Figure \ref{fig:M} shows two stereographic projections of permutahedron of $ Sp_6 $, where we have labelled each facet $ \gamma $ with the corresponding values of $ M_\gamma$ and $M'_\gamma$.

To see if $ P(M'_\bigdot) $ is an MV polytope, we check the tropical Pl\"ucker relations.  In particular, there is the tropical Pl\"ucker relation (see (\ref{eq:A2trop}))
\begin{equation*}
M'_3 + M'_{1-2} = \min( M'_{-2} + M'_{13} , M'_{1} + M'_{-23}).
\end{equation*}
However the left hand side of this equation equals $-2 $ and the right hand side equals $ -3 $ and so this equation does not hold.

Hence $ AM_1 \cdot P $ is not an MV polytope.

In fact, the actual $ f_1 \cdot P $ has BZ datum $ N_\bigdot $ with $ N_\gamma = M'_\gamma $ for all $ \gamma $ except that $ N_{1-2} = -3 $.  The facet corresponding to the chamber weight $ 1-2$ is shaded in figure \ref{fig:M}.  The vertices of $ e_1 \cdot P $ are $(0,0,0)$, $(0,2,0)$, $(-1,1,0)$, $(-1,2,1)$, $(0,0,x)$, $(-1,1,x)$, $(-1,2,x-1)$, $(0,2,x) $.

\subsection{Open questions}
This counterexample to the Anderson-Mirkovi\'c conjecture raises a number of interesting open problems.  

\begin{Question}
For which $(G, j) $ does the Anderson-Mirkovi\'c conjecture hold? 
\end{Question}

In the counterexample, we saw that $ AM_1 \cdot P \subset e_1 \cdot P $ (equivalently $ N_\gamma \le M'_\gamma $ for all $ \gamma $).  In fact, it is possible to prove that for all $ Sp_6 $ MV polytopes $P$, $ AM_1 \cdot P \subset e_1 \cdot P $.  This raises the following.

\begin{Question}
Is $ AM_j \cdot P $ always contained in $ f_j \cdot P $?  In other words, does $ f_j \cdot P $ satisfy the condition (iv) defining $AM_j $?
\end{Question}

One advantage of the Anderson-Mirkovi\'c conjecture is that via Proposition \ref{th:AMBZ} it gives a simple expression for $ M'_\bigdot $ (rather than recursively solving a sequence of $ (\min, +) $ equations).  So we ask:

\begin{Question}
Is there a simple closed expression for $ M'_\bigdot $ that holds in general?
\end{Question}

\section{Additional combinatorial structure}
So far in this paper we have constructed a crystal structure on the set $ \mathcal{P} $ of MV polytopes and proved that this crystal is isomorphic to the Verma crystal $ B(\infty) $.  The goal of this section is to recall some additional structure of $ B(\infty) $ and to explain how this structure can be expressed in terms of MV polytopes.

\subsection{Kashiwara involution}
Following Kashiwara \cite{Kas}, recall that there exists a weight preserving involution $ * : B(\infty) \rightarrow B(\infty) $.  This map corresponds to the anti-automorphism of $U^\vee_+ $ which fixes the Chevalley generators $ E_i $.  This involution gives rise to another crystal structure on the set $ B(\infty) $, with crystal operators $ f_i^* := * \circ f_i \circ *, e_i := * \circ f_i \circ * $.  One importance of this twisted crystal structure is that it is part of the data necessary to formulate the uniqueness theorem of Kashiwara-Saito \cite{KS} for the crystal $ B(\infty) $.

Any model for $ B(\infty)$, such as $ \mathcal{P}, \mathcal{B}, \mathcal{M} $, also acquires an involution $*$ and the twisted crystal structure.  

The following result is due to Lusztig.  In this form it appears as Proposition 3.3(iii) in \cite{BZtpm}.

\begin{Proposition}[{\cite{Lbookbook}}] \label{th:canstar}
Suppose that $ b \in \mathcal{B} $ and $ \wi $ is a reduced word.  Let $(n_1, \dots, n_m)= \phi_\wi(b) $.   Then, \begin{equation*}
\phi_{(\eta(i_m), \dots, \eta(i_1))}(*(b)) = (n_m, \dots, n_1)
\end{equation*}
\end{Proposition}

Here $ \eta : I \rightarrow I $ denotes the Dynkin diagram automorphism such that $ - w_0 \cdot \alpha_i = \alpha_{\eta(i)}$.

From this proposition, we can see immediately how to express $ * $ in terms of MV polytopes.

\begin{Theorem} \label{th:kashinv}
Let $ P \in \mathcal{P} $ denote a stable MV polytope.  Then $ *(P) = -P $.
\end{Theorem}

\begin{proof}
Suppose that $ P $ has BZ datum $ M_\bigdot$.  Note that $ -P $ is a pseudo-Weyl polytope defined by the collection $M_\bigdot' $ where $ M'_\gamma = M_{-\gamma} $ for all chamber weights $ \gamma $.  Checking the tropical Plucker relations and the non-degeneracy inequalities shows that $ -P $ is an MV polytope.  

Recall that the $\wi$-Lusztig datum of an MV polytope is defined by taking the lengths of the edges along a path through its 1-skeleton corresponding to $ \wi $.  When a polytope is negated the beginning of this path becomes its end and the directions of the edges are negated.  Hence it is now the path corresponding to the reduced word $ (\eta(i_m), \dots, \eta(i_1)) $.   So the theorem follows from Theorem \ref{th:can} and Proposition \ref{th:canstar}.
\end{proof}

Now that we can express the Kashiwara involution in terms of MV polytopes, it is easy to see how the crystal operators $ f_j^*, e_j^*$ act on MV polytopes.  In particular, Theorem \ref{th:can} and the above Theorem immediately imply the following result.

\begin{Corollary}
Let $ P $ be an MV polytope with GGMS datum $ \mu_\bigdot $.
\begin{enumerate}
\item $ f_j \cdot P $is the unique MV polytope whose set of vertices $ \mu'_\bigdot $ satisfies
\begin{equation*}
\mu'_{w_0} = \mu_{w_0} + \alpha_j^\vee \ \text{ and } \ \mu'_w = \mu_w \text{ if } s_j w > w.
\end{equation*}
\item $ e_j \cdot P = 0 $ if and only if $ \mu_{w_0 s_j} = \mu_{w_0} $.  Otherwise, $ e_j \cdot P $ is the unique MV polytope whose set of vertices $ \mu'_\bigdot $ satisfies
\begin{equation*}
\mu'_{w_0} = \mu_{w_0} - \alpha_j^\vee \ \text{ and } \ \mu'_w = \mu_w \text{ if } s_j w > w.
\end{equation*}
\end{enumerate}
\end{Corollary}
In other words, the $ f_j^* $ works at the highest vertex $ w_0 $ while $f_j $ works at the lowest vertex. 

\subsection{Finite crystals}
Recall that for any $ \lambda \in \Lambda_+ $, there exists a crystal $ B(\lambda) $ which is the crystal for the finite dimensional irreducible representation $V_\lambda $.  This crystal can also be described by MV polytopes.  
Let $ \mathcal{P}(\lambda) $ be the set of MV polytopes $ P = P(M_\bigdot) $ with highest vertex $ \lambda $ which satisfy one of the following two equivalent conditions:
\begin{enumerate}
\item $ M_{w_0 s_i \cdot \fund_i} \ge \langle w_0 \cdot \lambda, \fund_i \rangle $ for all $ i \in I $.
\item $ P \subset \conv(W \cdot \lambda) $.
\end{enumerate}
The equivalence of these two conditions was shown in \cite[section 8]{me}.

As shown in \cite[section 8]{me}, this is precisely the set of MV polytopes which indexes the MV basis for $V_\lambda $.  It is also the set of MV polytopes which indexes the canonical basis for $ V_\lambda $.  

Note that the set of stable MV polytopes is naturally in bijection with the set of MV polytopes with highest vertex $ \lambda $, so that we have an inclusion $ \mathcal{P}(\lambda) \hookrightarrow \mathcal{P} $. Define a crystal structure on $ \mathcal{P}(\lambda) $ by the following rules.  If $ P \in \mathcal{P}(\lambda) $ has coweight $ (\mu, \lambda) $, then we define $ \wt(P) := \mu $.  The crystal operators on $ \mathcal{P}(\lambda) $ are inherited from their action on $ \mathcal{P} $, except that we declare $ f_j \cdot P = 0 $ in $ \mathcal{P}(\lambda) $ if $ f_j \cdot P \notin \mathcal{P}(\lambda) $. 

From the above remarks and the general theory of how the crystal $ B(\lambda) $ sits inside the crystal $ B(\infty) $, the following result is immediate.

\begin{Theorem}
With this crystal structure, $ \mathcal{P}(\lambda) $ is isomorphic to the crystal $ B(\lambda) $.
\end{Theorem}

\subsection{Kashiwara data}
We now we recall the notion of Kashiwara (or string) data, first studied by Berenstein-Zelevinsky  \cite{BZstring} for $ \mathfrak{sl}_n $ and by Kashiwara \cite{Kas} for the general case.  Fix a reduced word $ \wi $ for $ w_0 $.  Let $ b \in B(\lambda) $ and define a sequence of non-negative integers $(p_1, \dots, p_m) $ by
\begin{align*}
p_1 &:= \max \{ p : f_{i_1}^p (b) \ne 0 \}, \\
p_2 &:= \max \{ p : f_{i_2}^p f_{i_1}^{p_1}(b) \ne 0 \}, \\
&\dots \\
p_m &= \max \{ p : f_{i_m}^p \cdots f_{i_1}^{p_1}(b) \ne 0 \}.
\end{align*}
In other words we apply lowering operators in the direction $ i_1 $ as far as we can, then in the direction $ i_2 $, etc.

The following result seems to be due to Littelmann \cite{Lcones}.
\begin{Proposition}
After we apply all these steps, we reach the lowest weight element of the crystal, i.e. $ f_{i_m}^{p_m} \cdots f_{i_1}^{p_1}(b) = b_\lambda^{low} $. 

Moreover, the map $ B(\lambda) \rightarrow \N^m $ taking $ b \mapsto (p_1, \dots, p_m) $ is injective.
\end{Proposition}

We call $ (p_1, \dots, p_m) $ the $\wi$-Kashiwara datum of $ b $.  In many ways, it is analogous to the notion of Lusztig data.  As we have an isomorphism $ B(\lambda) \simeq \mathcal{P}(\lambda) $, we may speak of the $\wi$-Kashiwara datum of an MV polytope.  We will now explain how this Kashiwara data can be easily recovered from the MV polytope.

Let $ P = P(M_\bigdot) = \conv(\mu_\bigdot) $ be an MV polytope.  Let $ w, i $ be such that $ ws_i > w $.  Then $ \mu_w $ and $ \mu_{ws_i} $ are connected by an edge in $ P $ which goes in the direction $ w \cdot \alpha_i^\vee $.  Recall (section \ref{se:MVpoly}) that lengths edges were used to determine the Lusztig datum of the MV polytope.  Now, we are more concerned with the following quantity
\begin{equation*}
\frac{1}{2} \langle \mu_w + \mu_{w s_i}, w \cdot \alpha_i \rangle \end{equation*}
which is the inner product of the midpoint of the edge with the direction of the edge.  We call this the \textbf{height of the midpoint} of the edge $ (\mu_w, \mu_{ws_i}) $.  A simple computation shows that we can express this quantity in terms of the BZ datum as
\begin{equation*}
\frac{1}{2} \langle \mu_w + \mu_{w s_i}, w \cdot \alpha_i \rangle = M_{w \cdot \fund_i} - M_{w s_i \cdot \fund_i}.
\end{equation*}

\begin{Theorem} \label{th:kashdata}
Let $ P = \conv(\mu_\bigdot) \in \mathcal{P}(\lambda) $ and let $\wi $ be a reduced word for $w_0$. The $ \wi$-Kashiwara data $(p_1, \dots, p_m) $ of $ P $ is the sequence of heights of the midpoints of the edges $ (\mu_e, \mu_{s_{i_1}}), \dots, (\mu_{s_{i_1} \cdots s_{i_{m-1}}}, \mu_{s_{i_1} \cdots s_{i_m}})$. In other words,
\begin{equation*}
p_k = M_{w^\wi_{k-1} \cdot \fund_{i_k}} - M_{w^\wi_k \cdot \fund_{i_k}}.
\end{equation*}
\end{Theorem}
(Recall our notation $w^\wi_k = s_{i_1} \cdots s_{i_k}$.)

The main step in the proof of this theorem is to compare the change of Kashiwara data with the tropical Pl\"ucker relations. Two reduced words $ \wi, \wi' $ are said to be related by a $d$-\textbf{braid move} involving $ i, j$, starting at position $ k $, if 
\begin{align*}
\wi &= (\dots, i_k, i, j, i, \dots, i_{k+d+1}, \dots), \\
\wi' &= (\dots, i_k, j, i, j, \dots, i_{k+d+1}, \dots), 
\end{align*}
where $ d $ is the order of $ s_i s_j $.

The $\wi$ and $\wi'$-Kashiwara data are related by the following Proposition, due to Berenstein-Zelevinsky \cite{BZstring} for $ d = 3$ and to Littelmann \cite{Lcones} and Nakashima \cite{N} for $ d = 4 $.

\begin{Proposition}\label{th:paramtrans}
Let $ \wi, \wi' $ be as above.
Suppose that $ (p_1, \dots, p_m), (p'_1, \dots, p'_m) $ are the $ \wi, \wi' $-Kashiwara data.

For $ l \notin \{k+1, \dots, k +d\} $, $ p_l = p'_l $.  For other $ l $ we have the following case by case formulas. 
\begin{enumerate}
\item If $ a_{ij} = 0 $, then $ d =2 $ and 
\begin{equation*}
p'_{k+1} = p_{k+2}, \ p'_{k+2} = p_{k+1}.
\end{equation*}

\item  If $ a_{ij} = a_{ji} = -1$, then $ d = 3 $ and
\begin{equation} \label{eq:paramtrans2}
\begin{gathered}
p'_{k+1} = \max (p_{k+3},p_{k+2} - p_{k+1}), \ p'_{k+2} = p_{k+1} + p_{k+3}, \\ p'_{k+3} = \min(p_{k+1}, p_{k+2} - p_{k+3}). \\
\end{gathered}
\end{equation}

\item If $a_{ij} = -1, a_{ji} = -2$, then $ d = 4$ and
\begin{equation} \label{eq:paramtrans3}
\begin{gathered}
p'_{k+1} = \max(p_{k+4}, p_{k+3} - p_{k+2}, p_{k+2} - p_{k+1}), \\ 
p'_{k+2} = \max(p_{k+3}, p_{k+1} - 2p_{k+2} + 2p_{k+3}, p_{k+1} + 2p_{k+4}), \\
p'_{k+3} = \min(p_{k+2},2p_{k+2} - p_{k+3} + p_{k+4},p_{k+4} + p_{k+1}), \\ 
p'_{k+4} = \min(p_{k+1}, 2p_{k+2} - p_{k+3}, p_{k+3} - 2p_{k+4}).
\end{gathered}
\end{equation}

\item If $ a_{ij} = -2, a_{ji} = -1$, then $ d = 4 $ and
\begin{equation} \label{eq:paramtrans4}
\begin{gathered}
p'_{k+1} = \max(p_{k+4}, 2p_{k+3} - p_{k+2}, p_{k+2} - 2p_{k+1}), \\
p'_{k+2} = \max(p_{k+3}, p_{k+1} - p_{k+2} + 2p_{k+3}, p_{k+1} + p_{k+4}), \\
p'_{k+3} = \min(p_{k+2},2p_{k+2} - 2p_{k+3} + p_{k+4},p_{k+4} + 2p_{k+1}), \\
p'_{k+4} = \min(p_{k+1}, p_{k+2} - p_{k+3}, p_{k+3} - p_{k+4}).
\end{gathered}
\end{equation}

\end{enumerate}

\end{Proposition}

Now, we compare this with MV polytopes.  This lemma is essentially due to Berenstein-Zelevinsky, \cite{BZtpm}, Theorem 5.2.
\begin{Lemma}
Let $ \wi, \wi', k$ be as above.
Suppose that $ M_\bigdot $ satisfies the tropical Pl\"ucker relation at $ (w_k^\wi, i, j) $.  Define $ (p_1, \dots, p_m) $ and $ (p'_1, \dots, p'_m) $ by the formulae
\begin{equation*}
p_l = M_{w^\wi_{l-1} \cdot \fund_{i_l}} - M_{w^\wi_l \cdot \fund_{i_l}}, \quad
p'_l = M_{w^{\wi'}_{l-1} \cdot \fund_{i'_l}} - M_{w^{\wi'}_l \cdot \fund_{i'_l}}
\end{equation*}
Then $(p_1, \dots, p_m) $ and $ (p'_1, \dots, p'_m) $ satisfy the conclusions of Proposition \ref{th:paramtrans}.
\end{Lemma}

\begin{proof}
For $ k < l \le k+d $, we must examine different cases depending on $a_{ij}, a_{ji} $.  We prove only the cases $ a_{ij} = 0 $ and $ a_{ij} = a_{ji} = -1 $.
\begin{enumerate}
\item If $ a_{ij} = 0 $, then $ d =2 $ and $ w^{\wi'}_{ k+2} = w^{\wi}_{ k+1}, w^{\wi'}_{ k+1} = w^{\wi}_{ k+2} $.  Hence it immediately follows that
\begin{equation*}
p'_{k+1} = p_{k+2}, \ p'_{k+2} = p_{k+1}.
\end{equation*}

\item  If $ a_{ij} = a_{ji} = -1$, then $ d = 3 $.  Let $ w = w_k^\wi $.  We have the tropical Pl\"ucker relation
\begin{equation*}
M_{w s_j \cdot \fund_j}  = \min \big( M_{w \cdot \fund_i} + M_{w s_i s_j \cdot \fund_j}, M_{w \cdot \fund_j} +  M_{w s_j s_i \cdot \fund_i} \big) - M_{w s_i \cdot \fund_i}.
\end{equation*}
Hence, we see that
\begin{equation*}
\begin{aligned}
 p'_{k+1} &= M_{w \cdot \fund_j} - M_{w s_j \cdot \fund_j} \\
 &= - \min \big( M_{w \cdot \fund_j} -M_{w \cdot \fund_i} - M_{w s_i s_j \cdot \fund_j} + M_{w s_1 \cdot \fund_i},  - M_{w s_j s_i \cdot \fund_i} + M_{w s_i \cdot \fund_i} \big) \\
&= \max(p_{k+2} - p_{k+1}, p_{k+3}) \\
\end{aligned}
\end{equation*}
as desired.
The formulae for $ p'_{k+2}, p'_{k+3} $ follow similarly.
\end{enumerate}
\end{proof}

\begin{proof}[Proof of Theorem \ref{th:kashdata}]
We proceed by induction starting with the lowest element $ b_\lambda^{low} $ as our base for the induction.

Note that $ b_\lambda^{low} $ is represented by the MV polytope $ \conv(W \cdot \lambda) $ which has BZ data $ M_{w \cdot \fund_i} = \langle \lambda, - w_0 \cdot \fund_i \rangle $ for all $ w, i$.  This is the unique MV polytope of weight $ (w_0 \cdot \lambda, \lambda) $ that fits inside $ \conv(W \cdot \lambda) $.  The heights of the midpoints of this polytope are all zero.  Likewise the Kashiwara data of $ b_\lambda^{low}$ are also all zero.

So now let $ P \in \mathcal{P}(\lambda) \setminus b_\lambda^{low} $ and suppose that the heights of the midpoints of the edges do give the Kashiwara data for all MV polytopes $ Q$ which are lower than $ P $ in the crystal.  By the previous lemma, it is enough to prove that the heights of the midpoints of the edges in $ P $ give the $ \wi$-Kashiwara datum of $ P $ for just one reduced word $ \wi $ (as we may pass from one reduced word to any other by a sequence of braid moves).  

Since $ P \ne b_\lambda^{low} $, there exists $ j \in I $ such that $ f_j (P) \ne 0 $.  Choose a reduced word $ \wi $ such that $ i_1 = j $.  Let $ (p'_1, \dots, p'_m) $ be the $\wi$-Kashiwara datum of $ Q := f_j (P) $.  From the definition of Kashiwara datum, it is immediate that  $(p_1, \dots, p_m) := (p'_1+1, \dots, p'_m) $ is the $\wi$-Kashiwara datum of $ P$.  Let $ M_\bigdot, M'_\bigdot$ be the BZ data of $ P, Q $ respectively.  By the induction hypothesis, 
\begin{equation*}
p'_k = M'_{w_{k-1}^\wi \cdot \fund_{i_k}} - M'_{w_{k}^\wi \cdot \fund_{i_k}}.
\end{equation*}

From the definition of the crystal structure on MV polytopes, we have that $ M_{\fund_i} = M'_{\fund_i} + 1 $ and $ M_\gamma = M'_\gamma $ for all other $ \wi$-chamber weights $ \gamma$.  Hence we see that 
\begin{equation*}
\begin{gathered}
p_1 = p'_1 + 1 = M'_{\fund_j} - M'_{s_j \cdot \fund_j} + 1 =  M_{\fund_j} - M_{s_j \cdot \fund_j}, \text{ and } \\
p_k = p'_k =  M'_{w_{k-1}^\wi \cdot \fund_{i_k}} - M'_{w_{k}^\wi \cdot \fund_{i_k}} = M_{w_{k-1}^\wi \cdot \fund_{i_k}} - M_{w_{k}^\wi \cdot \fund_{i_k}}, \text{ for } k > 1.
\end{gathered}
\end{equation*}
Therefore the result holds for $ P $ and so by induction the theorem is proved.
\end{proof}


\begin{thebibliography}{E-G-S}
\bibitem[A]{jared2}
J. E. Anderson, A polytope calculus for semisimple groups, \textit{Duke Math. J.} \textbf{116} (2003), no. 3, 567–-588; math.AG/0110225.

\bibitem[BZ1]{BZstring}
A. Berenstein and A. Zelevinsky, String bases for quantum groups of type $A_r$, \textit{Adv. in Soviet Math.} \textbf{16}, Part 1 (1993), 51--89.

\bibitem[BZ2]{BZschub}
A. Berenstein and A. Zelevinsky, Total positivity in Schubert varieties, \textit{Comment. Math. Helv.} \textbf{72} (1997), 128--166.

\bibitem[BZ3]{BZtpm}
A. Berenstein and A. Zelevinsky, Tensor product multiplicities, canonical bases, and totally positive varieties, \textit{Invent. Math.} \textbf{143} (2001), no. 1, 77--128; math.RT/\-9912012.

\bibitem[BFG]{BFG}
A. Braverman, M. Finkelberg, and D. Gaitsgory, Uhlenbeck Spaces via Affine Lie algebras, in ``Unity of Mathematics. In Honor of the Ninetieth Birthday of I.M.Gelfand'', \textit{Progr. in Math.} \textbf{244}, 2005; math.AG/0301176.

\bibitem[BG]{BG}
A. Braverman and D. Gaitsgory, Crystals via the Affine Grassmannian, \textit{Duke Math. J.}, \textbf{107} (2001) no. 3, 561--575; math.AG/9909077.

\bibitem[K]{me}
J. Kamnitzer, Mirkovi\'c-Vilonen cycles and polytopes; math.AG/0501365.

\bibitem[Kas]{Kas}
M. Kashiwara, The crystal base and Littelmann's refined Demazure character formula, \textit{Duke Math. J.} \textbf{71} (1993), 839--858.

\bibitem[KS]{KS}
M. Kashiwara and Y. Saito, Geometric construction of crystal basis, \textit{Duke Math. J.} \textbf{89} (1997), no. 1 9--36; q-alg/9606009.

\bibitem[Lit]{Lcones}
P. Littelmann, Cones, crystals, and patterns, \textit{Transform. Groups} \textbf{3} (1998), no. 2, 145--179.

\bibitem[L1]{Lbook}
G. Lusztig, Canonical bases arising from quantized enveloping algebras, \textit{ J. Amer. Math. Soc.} \textbf{3} (1990) no. 2, 447--498. 

\bibitem[L2]{Lbookbook}
G. Lusztig, \textit{Introduction to quantum groups}, Birkh\"auser, 1993.

\bibitem[M]{MG}
S. Morier-Genoud, Rel\`evement g\'eom\'etrique de la base canonique et involution de Sch\"utzenberger, \textit{C. R. Math. Acad. Sci. Paris}, \textbf{337} (2003), 371--374; math.RT/0309424.

\bibitem[N]{N}
T. Nakashima,  Polyhedral realizations of crystal bases and braid-type isomorphisms, \textit{Contemp. Math.}, \textbf{248} (1999), 419--435; math.QA/9906100.
\end{thebibliography}
\end{document}